\documentclass[prodmode,acmtecs]{acmsmall}
\usepackage{url}
\usepackage{stackengine} 
\usepackage{amsmath}
\usepackage{amssymb}
\usepackage[hyphens,spaces,obeyspaces]{trl}
\usepackage{listings}
\usepackage{tikz}
\usepackage{pgflibraryshapes}
\usetikzlibrary{backgrounds}
\usetikzlibrary{arrows}
\newenvironment{changemargin}[2]{%
  \begin{list}{}{%
    \setlength{\topsep}{0pt}%
    \setlength{\leftmargin}{#1}%
    \setlength{\rightmargin}{#2}%
    \setlength{\listparindent}{\parindent}%
    \setlength{\itemindent}{\parindent}%
    \setlength{\parsep}{\parskip}%
  }%
  \item[]}{\end{list}}

\graphicspath{{.}}

\newcommand{\diag}{\textrm{diag}}

\newcommand\OurNotes[1]{}
\newcommand\OurNotesR[1]{}

\newcommand{\dt}{{\Delta t}}

\title{PETSc/TS: A Modern Scalable ODE/DAE Solver Library}
\author{Shrirang Abhyankar, Jed Brown, Emil M. Constantinescu, Debojyoti Ghosh, \\ Barry F. Smith and Hong Zhang}
 
\begin{abstract}
  High-quality ordinary differential equation
(ODE) solver libraries have a long history, going back to the 1970s. Over the past
several years we have implemented, on top of the PETSc linear and nonlinear solver package, a new general-purpose, extensive,
extensible library for solving ODEs and differential algebraic equations (DAEs).
Package includes support for both forward and adjoint sensitivities
that can be easily utilized by the TAO optimization package, which
is also part of PETSc. The ODE/DAE integrator library strives to be highly
scalable but also to deliver high efficiency for modest-sized
problems. The library includes explicit solvers, implicit solvers, and
a collection of implicit-explicit solvers, all with a common user interface
and runtime selection of solver types, adaptive error control, and
monitoring of solution progress. The library also offers enormous
flexibility in selection of nonlinear and linear solvers, including the
entire suite of PETSc iterative solvers, as well as several parallel
direct solvers.
\end{abstract}




\begin{document}
\setcounter{page}{1}

\begin{bottomstuff}
Author's address: abhyshr@anl.gov, Energy Systems Division, Argonne
National Laboratory, 9700 South Cass Avenue, Lemont, IL 60439-5844; jed@jedbrown.org, Department of Computer Science,
University of Colorado Boulder, 430 UCB, Boulder, CO 80309;
emconsta@mcs.anl.gov, bsmith@mcs.anl.gov, and hongzhang@anl.gov Mathematics and Computer
Science Division, Argonne National Laboratory, 9700 South Cass Avenue,
Lemont, IL 60439-4844; ghosh5@llnl.gov, Center for Applied
Scientific Computing, Lawrence Livermore National Laboratory, 7000
East Avenue, Livermore, CA 94550.
\end{bottomstuff}

\maketitle

\OurNotes{Jed: Should we say something about nonlinear multigrid?}

\OurNotes{Jed: What about mentioning the petsc4py interface?  I think my Oregonator
  example is a nice starting place for people interested in ODE.}

\OurNotes{Lois: Do you want at least 1 sample computational result (maybe
  from a prior paper from Emil) to show the benefits of being able to
  use all of these different solvers.}

\OurNotes{Lois: Include some sample output from apps (e.g., how about
  a side-by-side science image and convergence plot showing something
  nifty with the algorithms + explaining why TS is essential here)} 

\OurNotes{Lois: Take a look at the TOMS optional reproducibility info.
  If including some computational results, it would be good to also
  satisfy the reproducibility requirements.} 

\OurNotes{Hong: Need to discuss applications more.}

\section{Introduction}
 Sophisticated numerical algorithms for the
integration of ordinary differential equations (ODEs) and differential
algebraic equations (DAEs) have existed for 
well over one hundred years, while general-purpose software libraries
for their solution have existed for at least forty years. With 
changes in the applications simulated and in the computer hardware,
such libraries are constantly evolving. Perhaps the best-known
such libraries for both ODEs and 
DAEs are those originating at Lawrence Livermore National
Laboratory, including VODE, CVODE, DASSL, \cite{DASSL} and
SUNDIALS \cite{Sundials}. Other libraries with ODE/DAE support include
Trilinos \cite{Trilinos}, as well as commercial codes such as Mathwork's MATLAB \cite{MATLAB}
and the Numerical Algorithms Group \cite{NAG}.
The theory of numerical algorithms for ODEs
and DAEs is covered in several monographs, including
\cite{Hairer_B2008_I,Hairer_B2002_II,Ascher_B1998,Brenan_B1996,Butcher_B2008}. 

The Portable, Extensible Toolkit for Scientific computation (Version
2.0), PETSc, was developed to provide scalable high-quality mathematical
libraries for distributed-memory parallel computing. Over the past several years we have
implemented a general-purpose, extensive, extensible ODE and DAE
integrator package with local and global error control, support for computing sensitivities, and handling of
events (discontinuities in solutions or parameters). In this paper we describe the design, properties, and usage
of the ODE/DAE integrators in PETSc. 
In Section 2 we introduce the organization of PETSc, followed in Section 3 
by the PETSc timestepping application programming interface (API)
and in Section 4 by the time integration schemes in PETSc. In Section 5 we 
discuss adaptive timestepping and error control and in Section 6 
the computation of sensitivities. In Section 7
we explain how events are handled and the tools for monitoring and
visualizing solutions and the solution process.
In Section 8 we discuss how PETSc handles discontinuities and events.
We conclude in Section 9
with a brief
summary of two high-level interfaces for accessing the ODE/DAE
integrators: one for networks (for example the power grid) and one for
chemical reactions.

\section{Background: PETSc Linear and Nonlinear Solvers}
\OurNotes{Lois: 
    Maybe include the diagram that we made last summer in section 2 near the list of IS, Vec, Mat, etc.  Or even better, maybe do a version of this that shows the TS perspective. }

\OurNotes{Lois:     I like Fig 1.  Could you make a counterpart of Figure 1 for TS generally?  A picture is worth a thousand words.}

PETSc is a scalable, MPI-based, object-oriented numerical software library written in C
and fully usable from C, C++, Fortran, and Python. See \cite{efficient} for 
details on its design and the users manual \cite{petsc-user-ref} for how to use  PETSC effectively. PETSc has several fundamental classes
from which applications are composed.
\begin{itemize}
\item {\tt IS} -- index sets used to index into vectors and matrices, an abstraction of a list of integers
\item {\tt Vec} -- vectors, abstract elements of $ R^n$, used to contain the ODE solutions, function evaluations, and so forth
\item {\tt Mat} -- matrices, representations of linear operations including matrix-free formulations and sparse and dense storage formats
\item {\tt PC} -- preconditioners, single-step iterative solvers
  including domain decomposition and algebraic and geometric multigrid as
  well as direct solvers such as LU
\item {\tt KSP} -- Krylov subspace solvers, multistep iterative solvers
\item {\tt SNES} -- nonlinear solvers, including Newton's method, quasi-Newton methods, and nonlinear Krylov methods
\end{itemize}
In addition PETSc has a variety of helper classes that are useful
for implicit ODE solvers. These include {\tt MatColoring} and {\tt MatFDColoring}, which are
used to efficiently compute Jacobians via finite difference; see Section \ref{computejacobinsviadifferencing}. Moreover,
PETSc has an abstract class {\tt DM} that serves as an adapter
between meshes, discretizations, and other problem descriptors and the
algebraic and timestepper objects that are used to solve the discrete problem.

PETSc takes a minimalist approach to public APIs, attempting to keep
them as small as possible and with as little overlap in functionalities
as possible. In addition PETSc provides both a programmatic
interface to almost all functionalities within PETSc and a simple
string-based system, called the options database, that allows runtime
control of almost all functionality in PETSc.

Because PETSc is written in C, which does not have native syntax to
create classes and class hierarchies of object-oriented languages,
the classes are managed ``manually'' by the use of C structs
and function pointers. One feature of this approach, which can also be
obtained in object-oriented languages through the use of delegators,
is that any object can be changed at runtime, even multiple times
during the run to different derived classes of the same base
class. For example, a linear solver object that uses the GMRES
method can later be changed to use Bi-CG-stab by simply resetting the solver type, without the need for factory classes.
Many PETSc functions
have optional arguments; since C does not support function overloading,
one passes {\tt PETSC\_DEFAULT} (for optional scalar arguments) and
{\tt NULL} (for optional array arguments). In order to allow changing
basic type sizes at compile time PETSc has its own types: {\tt
  PetscReal}, which can represent half-, single-, double-, or quad-precision floating point; {\tt PetscScalar}, which represents complex
values when PETSc is built with complex numbers; and {\tt PetscInt},
which represents either 32- or 64-bit integers. Since PETSc is written in C,
we cannot utilize templates for this purpose as would be done in C++.

The usage of PETSc objects generally proceeds in the following order.
\begin{itemize}
\item {\tt XXXCreate(MPI\_Comm comm,XXX *xxx)} creates an object of type XXX,
  for example, KSPCreate(MPI\_comm,KSP*) creates a Krylov solver object.
\item {\tt XXXSetYYY(XXX xxx,...)} sets options to the object, via the functional interface.
\item {\tt XXXSetType(XXX xxx,const char* typename)} sets the specific subclass of the object,
  for example, ``gmres'' is a subclass of KSP.
\item {\tt XXXSetFromOptions(XXX xxx)} allows setting the type and options of the object from the options database.
\item {\tt XXXSetYYY(XXX xxx,...)} sets additional options.
\item {\tt XXXSetUp(XXX xxx)} fully instantiates the object so that it is ready to be used.
\item Use the object.
\item {\tt XXXDestroy(XXX *xxx)} frees all the space being used by the solver. PETSc uses reference counting to ensure that objects referenced by multiple other objects are not prematurely destroyed.  
\end{itemize}
We present full examples of this approach below for the ODE/DAE solvers.

\section{PETSc Timestepping Application Programming Interface}

The PETSc interface for solving time-dependent problems is organized around the following form of a DAE:
\[
        F(t,u,\dot{u}) = G(t,u), \quad u(t_0) = u_0.
\]
If the matrix $F_{\dot{u}}(t) = \partial F
/ \partial \dot{u}$ is nonsingular, then the equation is an ODE and can be
transformed to the standard explicit form ($ \dot{u} = Q(t,u)$). This transformation
may not lead to efficient algorithms, so often the transformation to explicit form should be avoided. For ODEs with nontrivial mass
matrices such as those that arise in the finite element method, the implicit/DAE interface can significantly
reduce the overhead to prepare the system for algebraic solvers
by having the user assemble the correctly shifted matrix.  
This interface is also useful for ODE systems, not just DAE
systems.

To solve an ODE or DAE, one uses the timestep context {\tt TS} created
with {\tt TSCreate(comm,\&ts)} and then sets options from the options database with {\tt
  TSSetFromOptions(ts)}. To define the ODE/DAE, the user needs to provide one or more functions (callbacks). The TS API for providing these functions consists of the following.

\begin{itemize}
\item Function $F(t,u,\dot{u})$ is provided, by the user, with
\lstset{language=C,numbers=left,
    stepnumber=5,
    showstringspaces=false,
    tabsize=2,
    breaklines=true,
    breakatwhitespace=true}
\begin{lstlisting}
TSSetIFunction(TS ts,Vec r, (*f)(TS,PetscReal,Vec,Vec,Vec,void*),void *fP);
\end{lstlisting}

  The vector {\tt r} is an optional location to store the residual.
  The arguments to the function {\tt f()} are the
  timestep context, current time, input state $u$, input time derivative
  $\dot{u}$, and the (optional) user-provided context {\tt
    *fP} that contains data needed by the application-provided call-back routines. When only $ G(t,u)$ is provided, TS automatically assumes that $F(t,u,\dot{u}) = \dot{u} $.

\item Function $ G(t,u)$, if it is nonzero, is provided, by the user, with
\begin{lstlisting}
TSSetRHSFunction(TS ts,Vec r,(*g)(TS,PetscReal,Vec,Vec,void*),void *gP);
\end{lstlisting}
Again the vector {\tt r} is an optional location to store the residual.
  The arguments to the function {\tt g()} are the
  timestep context, current time, input state $u$, and the (optional) user-provided context {\tt *gP}.
  
\item Jacobian $(shift)F_{\dot{u}} + F_u$ \\
  If using a fully implicit or semi-implicit (IMEX) method, one also must provide an appropriate (approximate) Jacobian matrix of $F()$ at the current solution $u^n$ using
\begin{lstlisting}
TSSetIJacobian(TS ts,Mat A,Mat B,
    (*j)(TS,PetscReal,Vec,Vec,PetscReal,Mat,Mat,void*),void *jP);
\end{lstlisting}
  The arguments of {\tt j()}
  are the timestep context; current time; input state $u$; input
  derivative $\dot{u}$; $shift$ (described below); matrix $A$ (which defines the Jacobian); matrix $B$, which is optionally different from $A$ (from which the preconditioner is constructed);
  and the (optional) user-provided context {\tt jP}.
  
  This form for the Jacobian arises because for all implicit and semi-implicit time integrators in PETSc the value of $\dot{u}^n$ is approximated in the ODE/DAE solver algorithms by $(shift) u^n + q(u^{n-1},...)$,
  where the method-specific function $ q(u^{n-1},...)$ depends only on previous iterations. Hence
\begin{eqnarray*}
    \frac{d F}{d u^n} &  = &  \frac{\partial F}{\partial \dot{u}^n}\frac{\partial \dot{u}^n}{\partial u^n} + \frac{\partial F}{\partial u^n}\\
                      & = &   (shift) F_{\dot{u}^n}(t^n,u^n,\dot{u}^n) + F_{u^n}(t^t,u^n,\dot{u}^n).
\end{eqnarray*}
For example, the backward Euler method has $ \dot{u}^n = (u^n-u^{n-1})/\dt$. With 
$   F(u^n) = M \dot{u}^n - f(t,u^n)$, one
one obtains the expected Jacobian
\begin{eqnarray*}
  \frac{d F}{d u^n} &  = & \frac{\partial (M\dot{u}^n) }{\partial \dot{u}^n}\frac{\partial \dot{u}^n}{\partial u^n} - \frac{\partial f}{\partial u^n}\\
  & = & \frac{1 }{\dt} M 
  - f_{u^n}(t,u^n).
\end{eqnarray*}
In this case the value of \textit{shift} is $1/\dt$.

\item Jacobian $G_u$ \\ If using a fully implicit method and the
  function $ G() $ is provided, one must also provide an appropriate
  (approximate) Jacobian matrix of $G()$ using
\begin{lstlisting}
TSSetRHSJacobian(TS ts,Mat A,Mat B,
  (*gj)(TS,PetscReal,Vec,Mat,Mat,void*),void *gjP);
\end{lstlisting}
  The arguments for the function {\tt gj()}
  are the timestep context, current time, input state $u$, matrix $A$, optional matrix $B$ from which the preconditioning is constructed,
  and the (optional) user-provided context {\tt gjP}.
\end{itemize}

Providing appropriate $ F() $ and $G() $ and their derivatives for the problem allows for
easy runtime switching between explicit, semi-implicit, and
fully implicit methods.

\label{computejacobinsviadifferencing}

Providing correctly coded Jacobians is often a major stumbling block
for users of ODE/DAE integration packages. PETSc provides three
useful tools to help users in this process:
\begin{itemize}
  \item application of Jacobians via {\em matrix-free} differencing approaches,
  \item explicit computation of Jacobians via matrix coloring and differencing, and
  \item automatic testing of user-provided Jacobian computations.
\end{itemize}
Finite-difference-based matrix-free application of Jacobians is
handled with a special PETSc matrix class that never forms the matrix
entries explicitly but merely provides matrix-vector products. For
most situations the user simply provides the option {\tt -snes\_mf},
which uses the PETSc provided matrix-free matrix class, and either no preconditioner or a user-provided preconditioner or {\tt
  -snes\_mf\_operator}, where a standard preconditioner is constructed from
some user-provided approximation to the Jacobian. Users who desire more control over the process can utilize
\begin{lstlisting}
MatCreateMFFD(MPI_Comm,PetscInt m, PetscInt n,PetscInt M,PetscInt N,Mat *J)
MatMFFDSetFunction(Mat J,(*f)(void*,Vec,Vec),void *ctx)
\end{lstlisting}
The arguments of {\tt MatCreateMFFD()} are the local and global
dimensions of the operator, while the arguments of {\tt
  MatMFFDSetFunction()} include the nonlinear function and optional
user-provided context data. A simpler alternative uses the nonlinear
function already provided to the nonlinear solver with
\begin{lstlisting}
TSGetSNES(TS ts,SNES *snes)  
MatCreateSNESMF(SNES snes,Mat *J)
\end{lstlisting}

An explicit matrix representation of the Jacobian via
matrix coloring may be constructed by using the option {\tt -snes\_fd\_color}.
The coloring can be provided in complementary ways, either by
providing the nonzero structure of the Jacobian (but not its numerical
values) and applying a variety of matrix coloring routines to compute
the coloring (this is done by creating a coloring object with {\tt
  MatColoringCreate()} and from this performing the coloring) or by
  providing the coloring based on the mesh structure and specific
  numerical discretization approach used (this is done by calling {\tt
    DMCreateColoring()}). Once the coloring is provided, the actual
  computation of the Jacobian entries involves the use of {\tt
    MatFDColoringCreate()} and {\tt MatFDColoringSetFunction()}, which
  plays a role similar to {\tt MatMFFDSetFunction()}. Both the
  matrix-free differencing and the explicit computation of the
  Jacobian entries support various options for selecting the
  differencing parameters. The explicit computation of Jacobian
  entries via differencing can be used to find the locations of
  Jacobian entries incorrectly provided by the user. In the simplest case
  this is handled via the option {\tt -snes\_type test
    -snes\_test\_display}. Other options include {\tt
    -snes\_compare\_coloring} and {\tt
    -snes\_compare\_coloring\_display}.

We now present a simple, but complete, example code demonstrating the use of
TS to solve a small set of ODEs: $\dot{u}_0 = - \kappa u_0 u_1$,
$\dot{u}_1= - \kappa u_0 u_1$ $\dot{u}_2 = \kappa u_0 u_1$,
$u^{0}=[1.0,0.7,0.0]^T$.

\begin{lstlisting}
/*   Defines the ODE passed to the ODE solver   */
IFunction(TS ts,PetscReal t,Vec U,Vec Udot,Vec F,AppCtx *ctx){
  PetscScalar       *f;
  const PetscScalar *u,*udot;

  /*  Allow access to the vector entries  */
  VecGetArrayRead(U,&u); VecGetArrayRead(Udot,&udot); VecGetArray(F,&f);
  f[0] = udot[0] + ctx->k*u[0]*u[1];
  f[1] = udot[1] + ctx->k*u[0]*u[1];
  f[2] = udot[2] - ctx->k*u[0]*u[1];
  VecRestoreArrayRead(U,&u); VecRestoreArrayRead(Udot,&udot);
  VecRestoreArray(F,&f);
}
/*     Defines the Jacobian of the ODE passed to the ODE solver.     */
IJacobian(TS ts,PetscReal t,Vec U,Vec Udot,PetscReal a,Mat A,Mat B,AppCtx *ctx){
  PetscInt          rowcol[] = {0,1,2};
  PetscScalar       J[3][3];
  const PetscScalar *u,*udot;

  VecGetArrayRead(U,&u); VecGetArrayRead(Udot,&udot);
  J[0][0] = a + ctx->k*u[1];   J[0][1] = ctx->k*u[0];       J[0][2] = 0.0;
  J[1][0] = ctx->k*u[1];       J[1][1] = a + ctx->k*u[0];   J[1][2] = 0.0;
  J[2][0] = -ctx->k*u[1];      J[2][1] = -ctx->k*u[0];      J[2][2] = a;
  MatSetValues(B,3,rowcol,3,rowcol,&J[0][0],INSERT_VALUES);
  VecRestoreArrayRead(U,&u); VecRestoreArrayRead(Udot,&udot);
  MatAssemblyBegin(A,MAT_FINAL_ASSEMBLY);
  MatAssemblyEnd(A,MAT_FINAL_ASSEMBLY);
}
/*   Defines the initial conditions (and the analytic solution)  */
Solution(TS ts,PetscReal t,Vec U,AppCtx *ctx){
  const PetscScalar *uinit;
  PetscScalar       *u,d0,q;

  VecGetArrayRead(ctx->initialsolution,&uinit); VecGetArray(U,&u);
  d0   = uinit[0] - uinit[1];
  if (d0 == 0.0) q = ctx->k*t;
  else q = (1.0 - PetscExpScalar(-ctx->k*t*d0))/d0;
  u[0] = uinit[0]/(1.0 + uinit[1]*q);
  u[1] = u[0] - d0;
  u[2] = uinit[1] + uinit[2] - u[1];
  VecRestoreArrayRead(ctx->initialsolution,&uinit); VecRestoreArray(U,&u);
}
/* Creates the TS object, sets functions, options, then solves the ODE */
int main(int argc,char **argv){
  TS             ts;            /* ODE integrator */
  Vec            U;             /* solution will be stored here */
  Mat            A;             /* Jacobian matrix */
  PetscInt       n = 3;
  AppCtx         ctx;

  PetscInitialize(&argc,&argv,(char*)0,help);
  /*   Create necessary matrix and vectors  */
  MatCreate(PETSC_COMM_WORLD,&A);
  MatSetSizes(A,n,n,PETSC_DETERMINE,PETSC_DETERMINE);
  MatSetFromOptions(A); MatSetUp(A);
  MatCreateVecs(A,&U,&ctx.initialsolution);
  /*   Set runtime option  */
  ctx.k = .9;
  PetscOptionsGetScalar(NULL,NULL,"-k",&ctx.k,NULL);
  /*  Create timestepping solver context  */
  TSCreate(PETSC_COMM_WORLD,&ts);
  TSSetProblemType(ts,TS_NONLINEAR);
  TSSetType(ts,TSROSW);
  TSSetIFunction(ts,NULL,(TSIFunction) IFunction,&ctx);
  TSSetIJacobian(ts,A,A,(TSIJacobian)IJacobian,&ctx);
  /*  Set initial conditions */
  Solution(ts,0,U,&ctx);
  /*   Set solver options  */
  TSSetTimeStep(ts,.001);
  TSSetMaxSteps(ts,1000);
  TSSetMaxTime(ts,20.0);
  TSSetExactFinalTime(ts,TS_EXACTFINALTIME_STEPOVER);
  TSSetFromOptions(ts);
  TSMonitorLGSetVariableNames(ts,names);

  TSSolve(ts,U);

  VecDestroy(&ctx.initialsolution); MatDestroy(&A); VecDestroy(&U);
  TSDestroy(&ts);
  PetscFinalize();
}
\end{lstlisting}  

We next present a simple example code demonstrating the use of
TS to solve a small set of stiff ODEs (a 3-variable oscillatory ODE system from chemical
reactions, problem OREGO in \cite{Hairer_B2002_II}) written in Python and using the
petsc4py \cite{dalcin2011parallel} binding:
\begin{eqnarray*}
  \dot{u}_0 & = & - 77.27 (u_1 + u_0(1-8.375 \times  10^{-6} u_0 -u_1), \\
  \dot{u}_1 & = & \frac{1}{77.27} (u_2 - (1 + u_0)u_1), \\
  \dot{u}_2 & = & 0.161 (u_0-u_2), \\
 u(t=0) & =  & [1.0,2.0,3.0]^T.
\end{eqnarray*}

\OurNotesR{Emil:
  If we keep Orego.py, I'll generate a work precision diagram from python.}


\lstset{language=PYTHON,numbers=left,
    stepnumber=5,
    showstringspaces=false,
    tabsize=2,
    breaklines=true,
    breakatwhitespace=true}
\begin{lstlisting}
import sys, petsc4py
from matplotlib import pylab
from matplotlib import rc
import numpy as np
petsc4py.init(sys.argv)

from petsc4py import PETSc

class Orego(object):
    n = 3
    comm = PETSc.COMM_SELF
    def evalSolution(self, t, x):
        x.setArray([1, 2, 3])
    def evalFunction(self, ts, t, x, xdot, f):
        f.setArray([xdot[0] - 77.27*(x[1] + x[0]*(1 - 8.375e-6*x[0] - x[1])),
                    xdot[1] - 1/77.27*(x[2] - (1 + x[0])*x[1]),
                    xdot[2] - 0.161*(x[0] - x[2])])
    def evalJacobian(self, ts, t, x, xdot, a, A, B):
        B[:,:] = [[a - 77.27*((1 - 8.375e-6*x[0] - x[1]) -8.375e-6*x[0]),
                   -77.27*(1 - x[0]), 0],
                  [1/77.27*x[1], a + 1/77.27*(1 + x[0]),   -1/77.27],
                  [-0.161, 0, a + 0.161]]
        B.assemble()
        if A != B: A.assemble()
        return True # same nonzero pattern

OptDB = PETSc.Options()
ode = Orego()

J = PETSc.Mat().createDense([ode.n, ode.n], comm=ode.comm)
J.setUp()
x = PETSc.Vec().createSeq(ode.n, comm=ode.comm); f = x.duplicate()

ts = PETSc.TS().create(comm=ode.comm)
ts.setType(ts.Type.ROSW) # use Rosenbrock-W method

ts.setIFunction(ode.evalFunction, f)
ts.setIJacobian(ode.evalJacobian, J)

history = []
def monitor(ts, i, t, x):
    xx = x[:].tolist()
    history.append((i, t, xx))
ts.setMonitor(monitor)

ts.setTime(0.0)
ts.setTimeStep(0.1)
ts.setMaxTime(360)
ts.setMaxSteps(2000)
ts.setExactFinalTime(PETSc.TS.ExactFinalTime.INTERPOLATE)
ts.setMaxSNESFailures(-1) # allow unlimited failures (step will be retried)

# Set a different tolerance on each variable.
vatol = x.duplicate(array=[1e-2, 1e-1, 1e-4])
ts.setTolerances(atol=vatol,rtol=1e-3) # adaptive controller attempts to match this tolerance

snes = ts.getSNES()             # Nonlinear solver
snes.setTolerances(max_it=10)   # Stop nonlinear solve after 10 iterations (TS will retry with shorter step)
ksp = snes.getKSP()             # Linear solver
ksp.setType(ksp.Type.PREONLY)   # No Krylov method
pc = ksp.getPC()                # Preconditioner
pc.setType(pc.Type.LU)          # Use a direct solve

ts.setFromOptions()             # Apply run-time options, e.g. -ts_adapt_monitor -ts_type arkimex -snes_converged_reason
ode.evalSolution(0.0, x)
ts.solve(x)

if OptDB.getBool('plot_history', True):
    ii = np.asarray([v[0] for v in history])
    tt = np.asarray([v[1] for v in history])
    xx = np.asarray([v[2] for v in history])

    rc('text', usetex=True)
    pylab.suptitle('Oregonator: TS \\texttt{%s}' % ts.getType())
    pylab.subplot(2,2,1)
    pylab.subplots_adjust(wspace=0.3)
    pylab.semilogy(ii[:-1], np.diff(tt), )
    pylab.xlabel('step number')
    pylab.ylabel('timestep')

    for i in range(0,3):
        pylab.subplot(2,2,i+2)
        pylab.semilogy(tt, xx[:,i], "rgb"[i])
        pylab.xlabel('time')
        pylab.ylabel('$x_%d$' % i)
    pylab.show()
\end{lstlisting}

\begin{figure}[h]
  \centering
  \hspace*{-1cm}
  \begin{tabular}{cc}
    \includegraphics[scale=0.5]{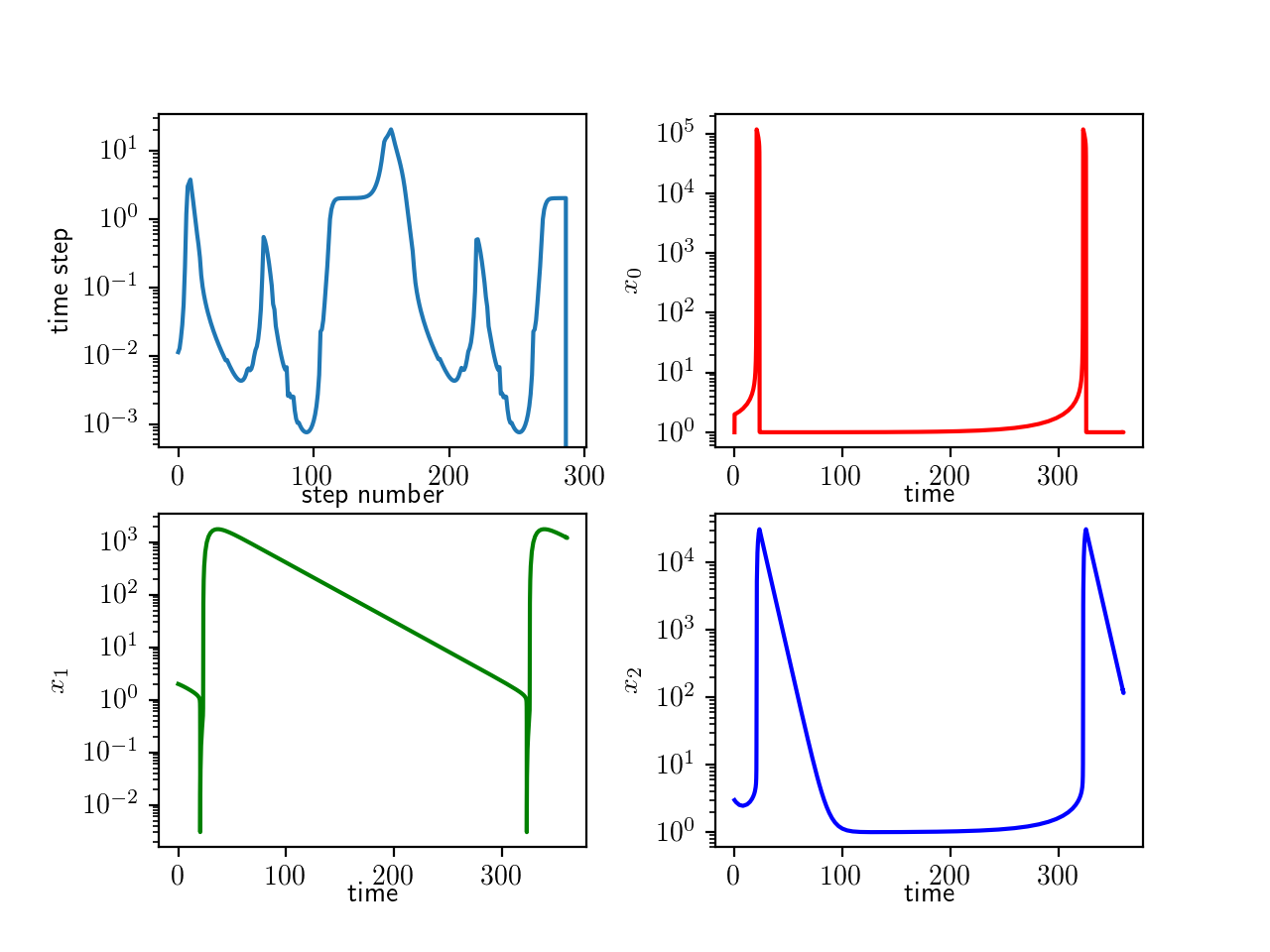}
    \hspace{-5mm}
  \includegraphics[scale=0.45]{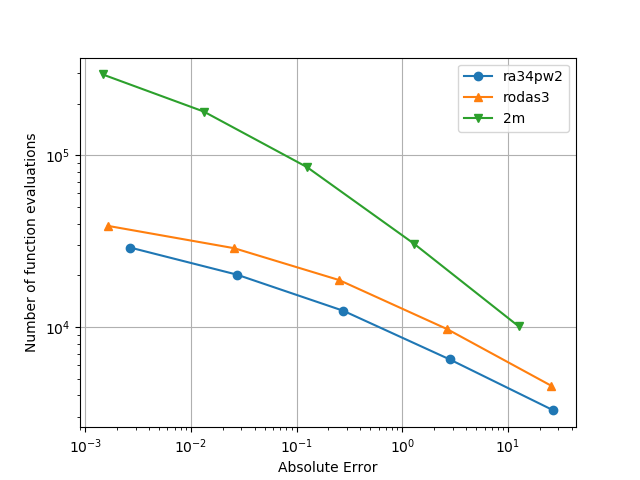}
  \end{tabular}
  \caption{Output and work-precision diagram for the OREGO problem with three
    Rosenbrock-W methods (\S\ref{sec:Ros}) under different TSAdapt
    tolerances}
  \label{fig:Orego}
\end{figure}

Figure \ref{fig:Orego} shows the output of the OREGO Python
code. We have also developed a work-precision diagram illustrating the
effect of choosing different tolerances in the TSAdapt (see
\S\ref{sec:Adaptors}) on the amount of effort and precision. 

\section{Timestepping Schemes}

 This section describes the interfaces
for setting the timestepping schemes and their options. The classes
of methods currently implemented in PETSc are described in Table
\ref{tab:TS:PETSc}. They cover multistage,  multistep, and general
linear methods with different stability properties. To address
different problem requirements, PETSc provides explicit methods that
are fast and accurate, implicit methods that have robust stability
properties, and partitioned methods that combine both implicit
and explicit integrators.
Many of
the methods implemented in PETSc allow users to register new schemes
by supplying a new set of method coefficients.  Most methods offer
local error control. Global error
estimation is also supported for several integrators.
When implicit and semi-implicit methods are used, any of the PETSc
linear and nonlinear solvers can be selected either by calling functions within the program or via the
PETSc options database. These are fully described in the PETSc users
manual \cite{petsc-user-ref}.
The following list details some of the
methods and their properties.
\begin{itemize}
\item {\bf euler} Explicit Euler
  method. This is a basic implementation of the simplest time integrator.
\item {\bf ssp} Class of strong-stability-preserving multistage explicit methods
  suitable for hyperbolic partial differential equations (PDEs). 
\item {\bf beuler, cn, theta}  low-order implicit
  methods suitable for DAEs and when
  stability is a concern.  
\item {\bf alpha(2)} Timestepping developed for Navier-Stokes
  problems \cite{Jansen_2000}. 
\item {\bf glle} Implementation of implicit general linear methods for stiff
  problems.   
\item {\bf rk} Implementation of explicit Runge-Kutta methods.
\item {\bf eimex} Implementation of extrapolated partitioned
  Runge-Kutta methods. These methods can have arbitrarily high
  orders \cite{Constantinescu_A2010a}.
\item {\bf arkimex} Implementation of additive partitioned Runge-Kutta
  methods. These are suitable for problems with stiff and nonstiff
  components. 
\item {\bf rosw} Implementation of Rosenbrock and W-methods,
  linearly implicit multistage methods with full or approximated
  Jacobian matrices. These methods are suitable for stiff, mildly
  nonlinear problems.
\item {\bf glee} Explicit and implicit general linear methods; typically
  self-starting, with global error estimation \cite{Constantinescu_2018a}. With suitable
  coefficients, these methods include euler, beuler, cn,
  theta, ssp, and rk. 
\item {\bf bdf} Standard backward differentiation methods. These are
  relatively low-order implicit multistep methods suitable for stiff
  problems and DAEs.
\end{itemize}



%

Explicit methods are conditionally stable. Implicit methods can be
conditionally or unconditionally stable. Unconditional stability can
be stronger or weaker. In this study we distinguish two types of
stability: \textit{A-Stable} methods, which have a stability region
that covers the entire real-negative complex half plane, and
\textit{L-Stable} or \textit{stiffly accurate} (SA) methods, which are
A-Stable methods for which the amplification factor goes to zero as
stiffness goes to infinity, thus giving them better stability
properties for stiff problems and DAEs.


%

\begin{table}[!htbp]
\centering
\tbl{\label{tab:TS:PETSc} List of time integration schemes
  available in PETSc.}
{  \small
\begin{tabular}{|c|p{4cm}|c|c|c|c|c|c|c|c|c|c|}
\hline
{\bf TS Name}&{\bf Reference}&{\bf Class}& {\bf Type}&{\bf Order}\\
\hline
euler&forward Euler&one-step&explicit&1\\
\hline
ssp&multistage SSP \cite{Ketcheson_2008}&Runge-Kutta&explicit&$\le$ 4\\
\hline
beuler&backward Euler&one-step&implicit&1\\
\hline
cn&Crank-Nicolson&one-step&implicit&2\\
\hline
theta&theta-method&one-step&implicit&$\le$2\\
\hline
alpha(2)&alpha-method \cite{Jansen_2000}&one-step&implicit&2\\
\hline
glle&general linear
\cite{Butcher_2007}&general linear&implicit& $\le$ 3\\
\hline
eimex&extrapolated IMEX \cite{Constantinescu_A2010a}&one-step&$\ge
1$, adaptive&\\
\hline
arkimex&\S \ref{sec:PRK} & IMEX Runge-Kutta & IMEX& $1-5$\\
\hline
rosw& \S \ref{sec:Ros} &Rosenbrock-W &linearly implicit&$1-4$\\
\hline
glee&  method with global error estimation \cite{Constantinescu_2018a} &general linear &explicit/implicit&$1-4$\\
\hline
bdf&  standard BDF methods \cite{Brenan_B1996} & multistep & implicit &
  $1-6$ \\
\hline
\end{tabular}}
\end{table}

\subsection{Partitioned Runge-Kutta\label{sec:PRK}}


Partitioned methods are aimed at relaxing ODE integrator classification into
strictly stiff and nonstiff problems that may have both types of
characteristics. Partitioned methods  employ two types of
integrators: an explicit integrator for the nonstiff problem
components and an implicit integrator suitable for the stiff
ones \cite{Ascher_1997,Kennedy_2003}. Partitioned Runge-Kutta methods
are typically strongly coupled; that is, both integrators participate in
each stage calculation. In the literature these methods are known as
implicit-explicit or semi-implicit \cite{Giraldo_2013,Zhong_1996}.  

A typical {\bf additive} partitioning of an ODE problem results in the
following:
\begin{align}
\label{eq:PETSC:ARK:partitioning}
\dot{u} =  \underset{\dot{U}_{E}}{\underbrace{G(t,u)}} +
\underset{\dot{U}_{I}}{\underbrace{H(t,u)}}\,, 
\end{align}
where $\dot{U}_{E}$ denotes the nonstiff right hand side function and $ F(t,u,\dot{u}) = \dot{u} - \dot{U}_{I}$
the stiff implicit function in PETSc. Integrating this problem explicitly in $G$ and
implicitly in $F$ from $u_n$ to $u_{n+1}$ by an additive Runge-Kutta
(ARK) method defined by coefficients $(A=\{a_{ij}\},b,c)$ for the
explicit part and $(\tilde{A},\tilde{b},\tilde{c})$ for the implicit
part, we have
\begin{subequations}
\begin{align}
U^{(i)} &= u_{n} + \dt \sum_{j=1}^{i-1} a_{ij} \dot{U}_{E}^{(j)} +
\dt \sum_{j=1}^{i} \tilde{a}_{ij} \dot{U}_{I}^{(j)}\,,~ i=1,\dots,s\\
\label{eq:ARK:ODE:compl}
u_{n+1} &= u_{n} + \dt \sum_{j=1}^{s} b_{j} \dot{U}_{E}^{(j)} +
\dt \sum_{j=1}^{s} \tilde{b}_{j} \dot{U}_{I}^{(j)}\,,
\end{align}
\end{subequations}
where $A=\{a_{ij}\}$ is strictly lower triangular,
$\tilde{A}=\{\tilde{a}_{ij}\}$ is lower triangular 
and can have zeros on the diagonal (these correspond to explicit
stages), and  $\circ^{(i)}$ is stage index $i$.

The implementation of the standard IMEX scheme is as follows. Solve
for $U^{(i)}$: 
\begin{subequations}
\label{eq:PETSC:ARKIMEX}
\begin{align}
\label{eq:PETSC:ARKIMEX:stage}
&F\left(t_n+c_{i}\dt,U^{(i)},\frac{1}{\dt\tilde{a}_{ii}}\left(U^{(i)}
-  Z\right)\right) =0\,,\\
\label{eq:PETSC:ARKIMEX:Z}
&  Z=u_{n} +
\dt \sum_{j=1}^{i-1} \tilde{a}_{ij} \dot{U}_{I}^{(j)} + 
 \dt \sum_{j=1}^{i-1} a_{ij} \dot{U}_{E}^{(j)}\\
\label{eq:PETSC:ARKIMEX:slopes}
&
\dot{U}_{I}^{(i)}=\frac{1}{\dt\tilde{a}_{ii}}\left(U^{(i)}
-  Z\right)\,,
\quad  \dot{U}_{E}^{(i)}=G\left(t_n+c_{i}\dt,U^{(i)}\right)\,, \quad
i=1,\dots,s\\
\label{eq:PETSC:ARKIMEX:completion}
& u_{n+1} = u_{n} + \dt \sum_{j=1}^{s} b_{j} \dot{U}_{E}^{(j)} +
\dt \sum_{j=1}^{s} \tilde{b}_{j} \dot{U}_{I}^{(j)}\,.
\end{align}
\end{subequations}
If $\tilde{a}_{ii}=0$, then \eqref{eq:PETSC:ARKIMEX:stage} is skipped,
and \eqref{eq:PETSC:ARKIMEX:slopes} is modified. This approach allows direct
use of these schemes for different types of problems as expressed in
Table \ref{tab_DE_forms}. 
%
%

Lower-order approximations are computed in the same way as for RK and ARK
methods by evaluating 
\eqref{eq:ARK:ODE:compl} with different $b$ and $\tilde{b}$ coefficients.

If one calls {\tt
  TSARKIMEXSetFullyImplicit()} or uses the option {\tt
  -ts\_imex\_fully\_implicit}, then \eqref{eq:PETSC:ARKIMEX:stage} solves
$F(t,u, \dot{u}) = G(t,u)$ by using only the implicit
integrator, thus turning the timestepping procedure into a
diagonally implicit integrator. This facilitates solving DAEs and
implicit ODEs. A summary of casting different problems through the
partitioned additive interface is given in Table  \ref{tab_DE_forms}.
An IMEX formulation for problems such as $M\dot{u}=g(t,u)+h(t,u)$
requires the user to provide $M^{-1} g(t,u)$.
 General cases such as
$F(t,u,\dot{u})=G(t,u)$ are not amenable to IMEX Runge-Kutta but can
be solved by using fully implicit methods, that is, by using the {\tt
  -ts\_imex\_fully\_implicit} option. 

\begin{table}[!htbp]
\centering
\tbl{\label{tab_DE_forms}Translation of various formulations of DAEs and ODEs into PETSc formulations.
 } 
    {\small
      \begin{tabular}{|l|p{1.3in}|p{2.4in}|}
        \hline
        $\dot{u}=g(t,u)$ & nonstiff ODE & $F(t,u,\dot{u}) := \dot{u}$,
        $G(t,u):=g(t,u)$\\
        \hline
        $\dot{u}=h(t,u)$ & stiff ODE & $F(t,u,\dot{u}) :=  \dot{u}-h(t,u)$,
        $G(t,u):={\tt NULL}$\\
        \hline
        $M \dot{u}=h(t,u)$ & stiff ODE with mass matrix& $F(t,u,\dot{u}) := M \dot{u}-h(t,u)$,
        $G(t,u):={\tt NULL}$\\
        \hline
        $M \dot{u} = g(t,u)$& nonstiff ODE with mass
        matrix&$F(t,u,\dot{u}) := \dot{u}, G(t,u) := M^{-1} g(t,u)$\\
        \hline
        $\dot{u}=g(t,u)+h(t,u)$ & stiff-nonstiff ODE & $F(t,u,\dot{u}) := \dot{u}-h(t,u)$, $G(t,u):=g(t,u)$\\
        $M\dot{u}=g(t,u)+h(t,u)$ & stiff-nonstiff ODE with mass matrix
        & $F(t,u,\dot{u}) := M \dot{u}-h(t,u)$,
        $G(t,u):=M^{-1}g(t,u)$\\
        \hline
        $h(t,y,\dot{y})=0$       & implicit ODE/DAE                  &
        $F(t,u,\dot{u}) := h(t,u,\dot{u})$, $G(t,u) := {\tt NULL}$;
            {\tt TSSetEquationType()} set to {\tt TS\_EQ\_IMPLICIT} or higher\\
            \hline
\end{tabular}}
\end{table}

The {\em dense output} or {\em continuous approximation} of the solution within
one timestep is also supported \cite{Horn_1983,Kennedy_2003}. This is used to obtain a high-order
interpolation on the one hand and a hot-start initial guess for the Newton
iterations on the other hand. 
The $s^* \ge s$ dense output formula for IMEX Runge-Kutta schemes of
order $p^*$ is given by  
\begin{align}
\label{eq:DenseOutpu:ARKt}
u^*(t_{n}+\theta \dt) := u_{n} + \dt \sum_{i=1}^{s^*}
B^*_i(\theta) g(t_{n} + c_i \dt, 
U^{(i)}) + \widehat{B}^*_i(\theta) f(t_{n} + c_i \dt,
U^{(i)})\,,
\end{align}
where  $\theta \in [0,1]$, $B^*_i(\theta)= \sum_{j=1}^{p^*} b^*_{ij}
\theta^j$, and $\widehat{B}^*_i(\theta)= \sum_{j=1}^{p^*} \widehat{b}^*_{ij}
\theta^j$.
We typically take $b^*_{ij}=\widehat{b}^*_{ij}$. When $\theta>1$, the
dense output is used for extrapolation. This option is set by {\tt
  -ts\_arkimex\_initial\_guess\_extrapolate} and has the effect of setting
the initial guess for all stages based on the dense output
extrapolated solution from the previous step. In nonlinear problems
this was shown to accelerate the code by up to three times; however,
the gains are highly problem dependent.  

\begin{table}[!htbp]
\tbl{\label{tab:IMEX:RK:PETSc} List of the IMEX RK schemes
  available in PETSc.}
    {  \stackengine{0pt}{
        \begin{tabular}{|c|c|c|c|c|c|c|c|c|}
\hline
{\bf TS}&{\bf Reference} &{\bf Stages} &{\bf Order}&{\bf Implicit}     &{\bf Stiff}\\
{\bf Name} &  &{\bf (IM)}   &{\bf (stage)}& {\bf stability} & {\bf accuracy}\\
\hline
1bee&B Euler +
Extrap&3(3)&1(1)&L&yes\\
\hline
a2&RK2a + Trap.&2(1)&2(2)&A&yes\\
\hline
l2&SSP2(2,2,2)\cite{Pareschi_2005}
&2(2)&2(1)&L&yes\\
\hline
ars122&ARS122,
\cite{Ascher_1997}&2(1)&3(1)&A&yes\\
\hline
2c&\cite{Giraldo_2013}&3(2)&2(2)&L&yes\\
\hline
2d&\cite{Giraldo_2013}&3(2)&2(2)&L&yes\\
\hline
2e&\cite{Giraldo_2013}&3(2)&2(2)&L&yes\\
\hline
prssp2&PRS(3,3,2)
\cite{Pareschi_2005}&3(3)&3(1)&L\\
\hline
3&\cite{Kennedy_2003}&4(3)&3(2)&L&yes\\
\hline
bpr3&\cite{Boscarino_TR2011}&5(4)&3(2)&L&yes\\
\hline
ars443&\cite{Ascher_1997}&5(4)&3(1)&L&yes\\
\hline
4&\cite{Kennedy_2003}&6(5)&4(2)&L&yes\\
\hline
5&\cite{Kennedy_2003}&8(7)&5(2)&L&yes\\
\hline
\end{tabular}
      }{
\begin{tabular}{|c|c|c|c|}
\hline
{\bf TS} &{\bf Embedded}&{\bf Dense}& {\bf Remarks}\\
{\bf Name} &    &{\bf  Output}& \\
\hline
1bee& yes(1)&no& extrapolated BEuler\\
\hline
a2&yes(1)&yes(2)&\\
\hline
l2&es(1)&yes(2)&SSP, SDIRK\\
\hline
ars122&yes(1)&yes(2)&\\
\hline
2c&yes(1)&yes(2)&SDIRK,SSP\\
\hline
2d&yes(1)&yes(2)&SDIRK\\
\hline
2e&yes(1)&yes(2)&SDIRK\\
\hline
prssp2&no&no&SSP, nonSDIRK\\
\hline
3&yes(2)&yes(2)&SDIRK \\
\hline
bpr3&no&no&SDIRK, DAE-1\\
\hline
ars443&no&no&SDIRK\\
\hline
4&yes(3)&yes(2,3)&SDIRK\\
\hline
5&yes(4)&yes(3)&SDIRK\\
\hline
\end{tabular}
      }{U}{c}{F}{F}{S}
   }
\end{table}

\subsection{Rosenbrock\label{sec:Ros}}
Rosenbrock methods are linearly implicit versions of
implicit Runge-Kutta methods. They use explicit function evaluations and implicit
linear solves, and therefore they tend to be faster than the implicit
Runge-Kutta methods because at each stage only a linear system needs
to be solved, 
as opposed
to the implicit Runge-Kutta methods that require solving a nonlinear
system at each stage. An $s$-stage Rosenbrock method is defined by 
coefficient matrices $A=a_{ij}\,,~j < i$ and
$\Gamma=\gamma_{i,j}\,,~j \le i$ and vector $b_i$, $i=1,\dots,s$. The
Rosenbrock scheme applied to $\dot{u}=f(t,u)$ computes the solution at step $n+1$ by
\begin{subequations}
\begin{align}
k_i&=\dt f(t_n + c_i \dt ,u_n+ \sum_{j=1}^{i-1}a_{ij} k_j) + \dt  J \sum_{j=1}^{i}
\gamma_{ij} k_j \,~ i=1,\dots,s\\
\label{eq:ros:comp}
u_{n+1} &=u_{n} +\sum_{i=1}^{s} b_i k_i\,,
\end{align}
\end{subequations}
where $J$ is the Jacobian matrix of $f(t,u)$ at $t=t_n$ and
$c_i=\sum_{j=1}^{i-1} a_{ij}$. Extensions to DAEs and PDAEs are readily
available \cite{Rang_2005}. The linear system is defined in terms
of the Jacobian matrix, which can be exact or approximated. The latter
case leads to W-methods.

\OurNotes{Jed: I think we should present Rosenbrocks in the contex of DAE, e.g., like
in Rang and Angermann (2007) Equation 2.1-2.1 (a simpler representation
and generalization of their 2005 paper, Equation 2.8-2.9).  Also,
looking at the code, it appears that we aren't actually computing with
this representation.  Is that true?  Why?} \OurNotesR{Emil: I agree,
  will do.}

We follow the implementation suggested by \cite{Rang_2005} and \cite{Kaps_1985}, where the
coefficient matrix $\Gamma$ is inverted and a change of variable is used:
\begin{align}
  \nonumber
v_i=\sum_{j=1}^{i} \gamma_{ij} k_j \,,~ i=1,\dots, s\,,
\end{align}
leading to the following
expressions:  
\begin{subequations}
  \label{eq:ros:impl}
\begin{align}
&\left(\frac{1}{\gamma_{ii} \dt} M - J\right) v_i = f\left(t_n + c_i
  \dt, v_n+ \sum_{j=1}^{i-1} \omega_{ij} v_j\right) + \frac{1}{\dt} M 
  \sum_{j=1}^{i-1} d_{ij} v_j
\,, ~i =1, \dots, s\\
&u_{n+1} =u_{n} +\sum_{j=1}^{s} m_j v_j\,,
\end{align}
\end{subequations}
where $\{d\}_{ij}=\diag(\gamma_{11}^{-1}, \dots,\gamma_{ss}^{-1}) - \Gamma
^{-1}$; $\{\omega\}_{ij}=A\Gamma^{-1}$; $\{m\}_i=b \,\Gamma^{-1}$; 
$\gamma_i=\sum_{j=1}^{i-1} \gamma_{ij}$; and $M$ is a mass matrix
that can be singular, resulting in a  DAE. In our implementation we
also allow for a noninvertible $\Gamma$ coefficient matrix by applying a
correction to \eqref{eq:ros:impl}. This allows us to use methods that
have explicit stages.
Lower-order approximations are computed in the same way as for RK and ARK
methods by evaluating 
\eqref{eq:ros:comp} with different $b$ coefficients. A work-precision
diagram with three of these methods is presented in
Fig. \ref{fig:Orego}.
%
\begin{table}
\centering
\tbl{\label{tab:IMEX:RosW:PETSc}List of the Rosenbrock W-schemes
  available in  PETSc.} 
     {  \stackengine{0pt}{
\begin{tabular}{|c|c|c|c|c|c|c|c|c|c|c|c|c|c|c|}
\hline
{\bf TS} &{\bf Reference}& {\bf Stages} &{\bf Order}&{\bf Implicit}     &{\bf Stiff} \\
{\bf Name} &  &{\bf (IM)}   &{\bf (stage)}&{\bf stability}& {\bf accuracy}\\
\hline
theta1& classical&
1(1)&1(1)&L&-\\
\hline
theta2& classical&
1(1)&2(2)&A&-\\
\hline
2m             &Zoltan&
2(2)&2(1)&L&No\\
\hline
2p             &Zoltan&
2(2)&2(1)&L&No\\
\hline
ra3pw          &\cite{Rang_2005}&
3(3)&3(1)&A&No\\
\hline
ra34pw2        &\cite{Rang_2005}&
4(4)&3(1)&L&Yes\\
\hline
rodas3         &\cite{Sandu_1997}&
4(4)&3(1)&L&Yes\\
\hline
sandu3         &\cite{Sandu_1997}&
3(3)&3(1)&L&Yes\\
\hline
assp3p3s1c &unpublished&
3(2)&3(1)&A&No\\
\hline
lassp3p4s2c&unpublished&
4(3)&3(1)&L&No\\
\hline
lassp3p4s2c&unpublished&
4(3)&3(1)&L&No\\
\hline
ark3       &unpublished&
4(3)&3(1)&L&No\\
\hline
\end{tabular}
    }{
\begin{tabular}{|c|c|c|c|c|c|c|c|c|c|c|c|c|c|c|}
\hline
{\bf TS} &{\bf Embedded}&{\bf Dense}&{\bf Inexact}&{\bf PDAE} &{\bf Remarks}\\
{\bf Name} &       &{\bf Output}& {\bf Jacobians}&&\\
\hline
theta1& &-&-&-&-\\
\hline
theta2&-&-&-&-&-\\
\hline
2m             &Yes(1)&Yes(2)&Yes&No&SSP\\
\hline
2p             &Yes(1)&Yes(2)&Yes&No&SSP\\
\hline
ra3pw          &Yes&Yes(2)&No&Yes(3)&-\\
\hline
ra34pw2        &Yes&Yes(3)&Yes&Yes(3)&-\\
\hline
rodas3         &Yes&No&No&Yes&-\\
\hline
sandu3         &Yes&Yes(2)&No&No&-\\ 
\hline
assp3p3s1c &Yes&Yes(2)&Yes&No&SSP\\
\hline
lassp3p4s2c&Yes&Yes(3)&Yes&No&SSP\\
\hline
lassp3p4s2c&Yes&Yes(3)&Yes&No&SSP\\
\hline
ark3       &Yes&Yes(3)&Yes&No&IMEX-RK\\
\hline
\end{tabular}
    }{U}{c}{F}{F}{S}
    }
\end{table}

For PDEs, much of the source code is
responsible for managing the mesh and spatial discretization, while
only a small amount handles the time integration. In PETSc the bridge
between the mass of code that handles the mesh and discretization and
the solver and time integrator is the {\tt DM} object. This object
provides the information needed by the solvers and integrators while
concealing all the details of the mesh and discretization
management. PETSc provides several DM classes including {\tt DMDA} for
structured grids with finite difference discretizations and {\tt
  DMPLEX} for unstructured meshes with finite element or finite volume
discretizations. We present an example of a PDE discretized by using finite differences on a two-dimensional
structured grid using the DM abstraction introduced earlier.
This example demonstrates an interesting nontrivial pattern formation
with a reaction-diffusion equation.
\OurNotes{Jed: The example on page 12-17 is much more verbose than necessary because it
uses the raw interfaces instead of interfaces such as
DMDATSSetRHSFunctionLocal.  I feel like the verbosity of the presently
chosen interfacse is sort of overwhelming for a journal article and may
not play well with reviewers.  Maybe this example can be reduced to just
show the local function and Jacobian, plus maybe an excerpt from
configuration.  Again, we should consider what comments are truly needed
for the audience reading this paper, rather than for someone learning
PETSc.}

\begin{lstlisting}
#include <petscdm.h>
#include <petscdmda.h>
#include <petscts.h>

typedef struct {
  PetscScalar u,v;
} Field;
typedef struct {
  PetscReal D1,D2,gamma,kappa;
} AppCtx;

int main(int argc,char **argv){
  TS             ts;                  /* ODE integrator */
  Vec            x;                   /* solution */
  DM             da;
  AppCtx         appctx;

  PetscInitialize(&argc,&argv,(char*)0,help);
  appctx.D1    = 8.0e-5;
  appctx.D2    = 4.0e-5;
  appctx.gamma = .024;
  appctx.kappa = .06;
  /*  Create distributed array (DMDA) to manage parallel grid and vectors */
  DMDACreate2d(PETSC_COMM_WORLD,DM_BOUNDARY_PERIODIC,DM_BOUNDARY_PERIODIC, DMDA_STENCIL_STAR,65,65,PETSC_DECIDE,PETSC_DECIDE,2,1,NULL,NULL,&da);
  DMSetFromOptions(da); DMSetUp(da);
  DMDASetFieldName(da,0,"u"); DMDASetFieldName(da,1,"v");
  DMCreateGlobalVector(da,&x);
  /* Create timestepping solver context */
  TSCreate(PETSC_COMM_WORLD,&ts);
  TSSetType(ts,TSARKIMEX);
  TSARKIMEXSetFullyImplicit(ts,PETSC_TRUE);
  TSSetDM(ts,da);
  TSSetProblemType(ts,TS_NONLINEAR);
  TSSetRHSFunction(ts,NULL,RHSFunction,&appctx);
  TSSetRHSJacobian(ts,NULL,NULL,RHSJacobian,&appctx);
  /* Set initial conditions */
  InitialConditions(da,x);
  TSSetSolution(ts,x);
  /*  Set solver options  */
  TSSetMaxTime(ts,2000.0);
  TSSetTimeStep(ts,.0001);
  TSSetExactFinalTime(ts,TS_EXACTFINALTIME_STEPOVER);
  TSSetFromOptions(ts);
  /*  Solve ODE system  */
  TSSolve(ts,x);
  VecDestroy(&x); TSDestroy(&ts); DMDestroy(&da);
  PetscFinalize();
}
/*   RHSFunction - Evaluates nonlinear function, F(x). */
RHSFunction(TS ts,PetscReal ftime,Vec U,Vec F,void *ptr) {
  AppCtx         *appctx = (AppCtx*)ptr;
  DM             da;
  PetscInt       i,j,Mx,My,xs,ys,xm,ym;
  PetscReal      hx,hy,sx,sy;
  PetscScalar    uc,uxx,uyy,vc,vxx,vyy;
  Field          **u,**f;
  Vec            localU;

  TSGetDM(ts,&da);
  DMGetLocalVector(da,&localU);
  DMDAGetInfo(da,PETSC_IGNORE,&Mx,&My,PETSC_IGNORE,...)
  hx = 2.50/(PetscReal)(Mx); sx = 1.0/(hx*hx);
  hy = 2.50/(PetscReal)(My); sy = 1.0/(hy*hy);
  /* Scatter ghost points to local vector  */
  DMGlobalToLocalBegin(da,U,INSERT_VALUES,localU);
  DMGlobalToLocalEnd(da,U,INSERT_VALUES,localU);
  DMDAVecGetArrayRead(da,localU,&u);
  DMDAVecGetArray(da,F,&f);
  /*  Get local grid boundaries
  DMDAGetCorners(da,&xs,&ys,NULL,&xm,&ym,NULL);
  /*  Compute function over the locally owned part of the grid */
  for (j=ys; j<ys+ym; j++) {
    for (i=xs; i<xs+xm; i++) {
      uc        = u[j][i].u;
      uxx       = (-2.0*uc + u[j][i-1].u + u[j][i+1].u)*sx;
      uyy       = (-2.0*uc + u[j-1][i].u + u[j+1][i].u)*sy;
      vc        = u[j][i].v;
      vxx       = (-2.0*vc + u[j][i-1].v + u[j][i+1].v)*sx;
      vyy       = (-2.0*vc + u[j-1][i].v + u[j+1][i].v)*sy;
      f[j][i].u = appctx->D1*(uxx + uyy)-uc*vc*vc+appctx->gamma*(1.0-uc);
      f[j][i].v = appctx->D2*(vxx + vyy)+uc*vc*vc-(appctx->gamma + appctx->kappa)*vc;
    }
  }
  DMDAVecRestoreArrayRead(da,localU,&u);
  DMDAVecRestoreArray(da,F,&f);
  DMRestoreLocalVector(da,&localU);
}
RHSJacobian(TS ts,PetscReal t,Vec U,Mat A,Mat BB,void *ctx) {
  AppCtx         *appctx = (AppCtx*)ctx;  /* application context */
  DM             da;
  PetscInt       i,j,Mx,My,xs,ys,xm,ym;
  PetscReal      hx,hy,sx,sy;
  PetscScalar    uc,vc;
  Field          **u;
  Vec            localU;
  MatStencil     stencil[6],rowstencil;
  PetscScalar    entries[6];

  TSGetDM(ts,&da);
  DMGetLocalVector(da,&localU);
  DMDAGetInfo(da,PETSC_IGNORE,&Mx,&My,PETSC_IGNORE,...)
  hx = 2.50/(PetscReal)(Mx); sx = 1.0/(hx*hx);
  hy = 2.50/(PetscReal)(My); sy = 1.0/(hy*hy);
  DMGlobalToLocalBegin(da,U,INSERT_VALUES,localU);
  DMGlobalToLocalEnd(da,U,INSERT_VALUES,localU);
  DMDAVecGetArrayRead(da,localU,&u);
  DMDAGetCorners(da,&xs,&ys,NULL,&xm,&ym,NULL);

  stencil[0].k = 0; stencil[1].k = 0; stencil[2].k = 0;
  stencil[3].k = 0; stencil[4].k = 0; stencil[5].k = 0;
  rowstencil.k = 0;
  for (j=ys; j<ys+ym; j++) {
    stencil[0].j = j-1;
    stencil[1].j = j+1;
    stencil[2].j = j;
    stencil[3].j = j;
    stencil[4].j = j;
    stencil[5].j = j;
    rowstencil.k = 0; rowstencil.j = j;
    for (i=xs; i<xs+xm; i++) {
      uc = u[j][i].u;
      vc = u[j][i].v;
      uyy       = (-2.0*uc + u[j-1][i].u + u[j+1][i].u)*sy;
      vxx       = (-2.0*vc + u[j][i-1].v + u[j][i+1].v)*sx;
      vyy       = (-2.0*vc + u[j-1][i].v + u[j+1][i].v)*sy;
      f[j][i].u = appctx->D1*(uxx + uyy)-uc*vc*vc+appctx->gamma*(1.0-uc);

      stencil[0].i = i; stencil[0].c = 0; entries[0] = appctx->D1*sy;
      stencil[1].i = i; stencil[1].c = 0; entries[1] = appctx->D1*sy;
      stencil[2].i = i-1; stencil[2].c = 0; entries[2] = appctx->D1*sx;
      stencil[3].i = i+1; stencil[3].c = 0; entries[3] = appctx->D1*sx;
      stencil[4].i = i; stencil[4].c = 0; entries[4] = -2.0*appctx->D1*(sx + sy) - vc*vc - appctx->gamma;
      stencil[5].i = i; stencil[5].c = 1; entries[5] = -2.0*uc*vc;
      rowstencil.i = i; rowstencil.c = 0;
      MatSetValuesStencil(A,1,&rowstencil,6,stencil,entries,INSERT_VALUES);

      stencil[0].c = 1; entries[0] = appctx->D2*sy;
      stencil[1].c = 1; entries[1] = appctx->D2*sy;
      stencil[2].c = 1; entries[2] = appctx->D2*sx;
      stencil[3].c = 1; entries[3] = appctx->D2*sx;
      stencil[4].c = 1; entries[4] = -2.0*appctx->D2*(sx + sy) + 2.0*uc*vc - appctx->gamma - appctx->kappa;
      stencil[5].c = 0; entries[5] = vc*vc;
      rowstencil.c = 1;
      MatSetValuesStencil(A,1,&rowstencil,6,stencil,entries,INSERT_VALUES);
    }
  }
  DMDAVecRestoreArrayRead(da,localU,&u);
  DMRestoreLocalVector(da,&localU);
  MatAssemblyBegin(A,MAT_FINAL_ASSEMBLY);
  MatAssemblyEnd(A,MAT_FINAL_ASSEMBLY);
  MatSetOption(A,MAT_NEW_NONZERO_LOCATION_ERR,PETSC_TRUE);
}
\end{lstlisting}
    
\subsection{Adaptive Timestepping and Error Control
 } \label{sec:Adaptors}

PETSc provides several options for automatic timestep control in order to
attain a user-specified goal via a {\tt TSAdapt} context. Typically,
the goals are related to accuracy. In this case the user 
provides an absolute (ATOL) and a relative (RTOL) error tolerance. The
adaptor controls the timestep in order to meet the specified error
tolerance. Most timestepping methods with adaptivity evaluate a lower-order
approximation at each timestep by using a different set of
coefficients Denote this solution as $\tilde{u}$. The following
weighted error quantity is used for timestep control:
\begin{align}
werr(t_{[n]}) = \frac
{||u(t_{[n]}) - \tilde{u}(t_{[n]})||_{\{1,2,\dots,\infty\}}}
{{\rm ATOL} + {\rm RTOL} \max(|u(t_{[n]})|,|\tilde{u}(t_{[n]})|)}\,.
\end{align}
If $werr(t_{[n]})$ is larger than one, then the estimated local
trunctation error at
step $n$ exceeds ATOL or RTOL. Otherwise, the estimated error is less
than that prescribed by the user, in which case 
the step is accepted and the next step adjusted so
that it tracks whether $werr(t_{[n]})$ will approach the value one. If the
error exceeds the tolerances 
specified by the user, then the step is rejected, and a smaller timestep is taken. 
This logic is implemented in the ``basic''
adaptor. A more advanced adaptivity logic based on linear digital
control theory and aimed at producing smoother step size sequences is
implemented in the ``dsp'' adapter \cite{Soderlind_2003,Soderlind_2006}.

In many fluid dynamics applications the timestep is restricted by
stability considerations as given by the Courant-–Friedrichs-–Lewy (CFL)
condition. TS provides an adapter that controls the timestep so that the
CFL stability is not exceeded.
Additionally, a special adapter for controlling the global error for the {\tt TS}
glee method \cite{Constantinescu_2018a} is available. This adapter can be used wherever the
standard (basic) one is used. Similar to the basic adapter, the
glee adapter can be used for tracking the absolute and relative
errors separately.

A list of timestep adapters is presented in Table
\ref{tab:TS:Adapt}. Custom adapters can be easily registered via the PETSc API.
\begin{table}[!htbp]
\centering
\tbl{\label{tab:TS:Adapt} List of time integration adapter schemes
  available in  PETSc.}
{  \scriptsize
\begin{tabular}{|p{2cm}|p{4cm}|p{4cm}|}
\hline
{\bf TS Adapt Name}&{\bf Remarks}&{\bf Used by}\\
\hline
none&No adaptivity & all\\
\hline
basic&Standard timestep adaptivity \cite{Gear_B1971} & all with lower-order
error approximation\\
\hline
dsp&Adapter using control theory \cite{Soderlind_2003} & same as basic\\
\hline
cfl&Controls the timestep to match error provided CFL limit&
typically {\tt TS} spp, rk\\
\hline 
glee&Time step adaptivity with global error
estimation \cite{Constantinescu_2018a}&typically for {\tt TS} glee methods, extends {\tt
  TSAdapt} basic\\ 
\hline
\end{tabular}}
\end{table}

\section{Computing Sensitivities (Derivatives)}

The timestepping library provides a framework based on discrete forward
(tangent linear) and adjoint models for sensitivity analysis for ODEs and DAEs.
The ODE/DAE solution process (henceforth called the forward run) can be obtained
by using either explicit or implicit solvers in TS, depending on the problem
properties. Currently supported method types are TSRK (Runge-Kutta) explicit
methods and TSTHETA (Theta) implicit methods.

\subsection{Discrete adjoint sensitivity}
\newcommand{\I}{\mathbf{\mathcal{I}}}
\newcommand{\N}{\mathbf{\mathcal{N}}}
\newcommand{\Sen}{\mathcal{S}}

The TSAdjoint routines of PETSc provide the capability to calculate partial derivatives of a given objective function
\begin{align}
  \Psi_i(u_0,p) = \Phi_i(u_F,p) + \int_{t_0}^{t_F} r_i(u(t),p,t)dt \quad i=1,...,n_\text{obj},
  \label{def:obj}
\end{align}
subject to initial conditions $u_0$ and parameters $p$.

Without loss of generality, we assume that the system is integrated with a one-step
method, 
\begin{equation}
  \label{eq:TSNumerical}
  u_{n+1} = \N_n(u_n), \quad n=0,\dots,N-1,\quad u_0=\I,
\end{equation}
where $\I$ are the initial values and the solution at the end of the simulation is given by $u_N$.

To illustrate the approach, we consider a simple
case in which we compute the sensitivities of the terminal function $\psi(u_N)$ to
initial values only. We use the Lagrange multipliers
$\lambda_0$,$\dots$,$\lambda_{N}$ to account for the constraint from each 
timestep, and we define the Lagrangian as
\OurNotes{Jed: The Lagrangian $\mathcal L$ in \eqref{eqn:lag} is not just a function of $I$, but also
$u$ and $\lambda$.  Then we can write (11) using partial derivatives and the
chain rule.  At least that is how I use Lagrange multipliers and I think
leads to clear presentation.}
\begin{equation}
\mathcal{L}(\I,u_0,\dots,u_N,\lambda_0,\dots,\lambda_N) = \psi(u_N)  - \lambda_0 ^T \left( u_0 - \I \right) - \sum_{n=0}^{N-1} \lambda_{n+1} ^T \left( u_{n+1} - \N (u_n) \right).
\label{eqn:lag}
\end{equation}

Differentiating Equation \eqref{eqn:lag} w.r.t. $\I$ and applying the chain rule, we obtain 
\begin{align}
\dfrac{d \mathcal{L} }{d \I} = \lambda_0 ^T   - \left( \frac{d \psi}{d u}(u_N) -  \lambda_{N} ^T \right)  \frac{\partial u_N}{\partial \I}
- \sum_{n=0}^{N-1}\left( \lambda_{n}^T - \lambda_{n+1} ^T \frac{d \N }{d y} (u_n) \right)  \frac{\partial u_n}{\partial \I}.
\end{align}
By defining $\lambda$ to be the solution of the discrete adjoint model,
\begin{align}
 \lambda_{N} = \left( \frac{d \psi}{d u}(u_n) \right)^T , \quad \lambda_{n}  =  \left( \frac{d \N}{d u} (u_n)  \right)^T \lambda_{n+1}, \quad n= N-1, \dots, 0,
\label{eqn:disadjoint}
\end{align} 
we obtain the gradient
$
\nabla_{\I} \mathcal{L} =  \nabla_{\I} \psi(u_n) =  \lambda_0 .
$

This model can be expanded to incorporate integral objective functions and calculate parametric sensitivities by augmenting the state variable into a larger system. See \cite{Zhang_2017ab} for more details.

To efficiently calculate the gradient with the adjoint method, one needs to first
perform a forward run that solves the original equation and saves the solution
trajectory with a checkpointing scheme, initialize the adjoint
sensitivity variables, and then perform a backward run that propagates the
adjoint sensitivity according to \eqref{eqn:disadjoint}. As can be seen from equation \eqref{eqn:disadjoint}, performing an adjoint step requires
trajectory information including the solution vector at the current step and
optional stage values if a multistage timestepping method is used. Applying
checkpointing techniques that have partial recomputation of the solution
provides a balance between recomputation and storage.

To use the PETSc adjoint solver, one creates two arrays of $n_\text{cost}$
vectors $\lambda$ and $\mu$ (if there are no parameters $p$, $\mu$ can be set to
NULL). The $\lambda$ vectors have the same dimension and parallel layout as the
solution vector for the ODE, and the $\mu$ vectors are of dimension $p$. The vectors $\lambda_i$
and $\mu_i$ should be initialized with the values $d\Phi_i/dy|_{t=t_F}$ and
$d\Phi_i/dp|_{t=t_F}$ respectively.

If $F()$ is a function of $p$, one needs to also provide the Jacobian $F_p$ with
\begin{lstlisting}
TSAdjointSetRHSJacobian(TS ts,Mat Amat, (*fp)(TS,PetscReal,Vec,Mat,void*),void *ctx)
\end{lstlisting}
The arguments for the function fp() are the timestep context, current time, $u$,
and the (optional) user-provided context.
If there is an integral term in the cost function, one must also provide
Jacobian terms for the integrand with
\begin{lstlisting}
TSSetCostIntegrand(TS ts,PetscInt numcost,Vec costintegral, (*rf)(TS,PetscReal,Vec,Vec,void*),(*drduf)(TS,PetscReal,Vec,Vec*,void*), (*drdpf)(TS,PetscReal,Vec,Vec*,void*),void *ctx)
\end{lstlisting}
where $\mathrm{drduf}= dr /du$, $\mathrm{drdpf} = dr /dp$.
The integral term can be evaluated in either the forward run or the backward
run by using the same timestepping algorithm as for the original equations.
\begin{figure}[ht]
  \centering
  \includegraphics[scale=0.4]{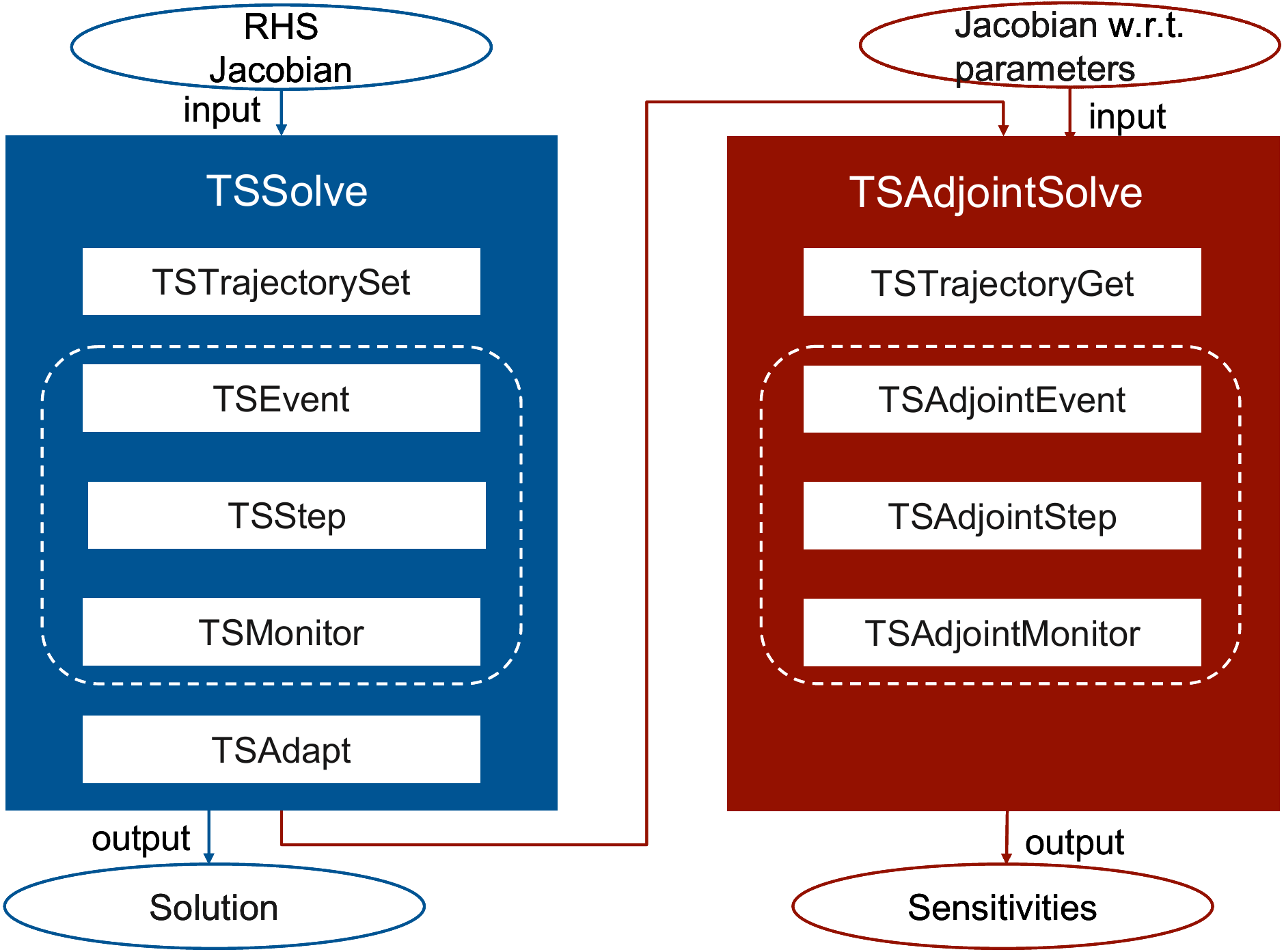}
  \caption{Structure of PETSc implementation for adjoint sensitivity analysis~\cite{Zhang_2017ab}}
  \label{fig:petsc_impl}
\end{figure}

The features of the PETSc adjoint solver are summarized as follows.
\begin{itemize}
  \item Supports multiobjective sensitivity calculation and integral objective functions
  \item Handles hybrid dynamical systems with discontinuities 
  \item Contains state-of-the-art checkpointing schemes
\end{itemize}

The adjoint solver is built on existing components in PETSc's timestepping
library \texttt{TS}, as shown in Fig. \ref{fig:petsc_impl}.
\begin{enumerate}
\item The \texttt{TSEvent} object, further introduced in Sec.
\ref{sec:tsevent}, supports detecting events and allows users to add a
post-event handler to modify the right-hand side function, reinitialize the DAE
system, and apply jump conditions for sensitivity analysis. It is particularly
important for the simulation of hybrid dynamical systems.

\item The \texttt{TSTrajectory} object provides a variety of sophisticated
online and offline checkpointing schemes that are suitable for single-level
storage media (for example, RAM) and multilevel storage media (RAM and external
disk/tape). Trajectory information is stored as checkpoints in the forward
run by repeatedly calling \texttt{TSTrajectorySet} at each timestep.
\texttt{TSTrajectoryGet} is responsible for obtaining the required trajectory
information before an adjoint step starts. It may extract the information
from the restored checkpoint directly or recompute from the checkpoint.
Recomputation typically happens when checkpoints are stored only at selective
timesteps because of limited storage capacity. \texttt{TSTrajectorySet} and
\texttt{TSTrajectoryGet} encapsulate the state-of-the-art checkpointing scheduler
\texttt{revolve} \cite{griewank2000algorithm} that can generate a guaranteed
optimal strategy.

\item \texttt{TSAdjointStep} corresponds to the adjoint version of
\texttt{TSStep}, which fulfills the timestepping operator $\N(u_n)$. Thus they
have similar underlying infrastructure, and their implementations differ from one
timestepping method to another. By design, the inputs for the adjoint solver are
either reused or modified from the original \texttt{TS} solver.
\end{enumerate}
All the components are compatible with one another and used together
within the highly composable solver in order to tackle the
difficulties of hybrid systems.  Details on using the infrastructure
discussed here for solving PDE-constrained optimization problems
utilizing the spectral element method can be found in
\cite{pdeconstrainedspectraladjoints}.

\subsection{Discrete forward (tangent linear) sensitivity}
\OurNotes{Jed: I think "\S 6.2 Discrete forward methods" should either include
"tangent" or "sensitivity".}
The discrete forward (also known as tangent linear) model for a one-step time
integration algorithm can be obtained by taking the derivative of
\eqref{eq:TSNumerical} with respect to the parameters.
The propagation equation for parameters $p$ can be symbolically described by
\begin{equation}
 \label{eqn:method1fwd}
   \Sen_0 = \frac{d \I}{d p}, \quad  \Sen_{n+1}  =  \frac{d \N}{d u} (u_n) \Sen_n, \quad n= 0, \dots, N-1,
\end{equation}
where $\Sen_{n} = d X_n /dp$ is a matrix denoting the solution
sensitivities (or so-called trajectory sensitivities in the power system field).
Note that each parameter results in one corresponding column of the sensitivity
matrix $\Sen$ and one linear equation to be solved at each timestep. 
Consequently, the computational cost of the forward approach is linear in the number
of parameters for which the sensitivities are calculated. This feature usually
limits its application to cases involving few parameters.

Like the discrete adjoint models, the implementation of discrete forward models
also depends on the particular time integration algorithm. In principle, these
two models are analogous to the well-known forward and reverse modes of
algorithmic differentiation (AD) that are applied to high-level abstractions of
a computer program. Traditional AD handles a sequence of operations (either a
source code line or a binary instruction) while in our case the primitive
operation is a timestep.

Furthermore, the forward model requires the same ingredients as those needed in
the adjoint model. Users may need to provide {\tt TSAdjointSetRHSJacobian()} and
{\tt TSSetCostIntegrand()} in the same way that they are used for TSAdjoint.

Although forward sensitivities are not used as frequently with gradient-based
optimization algorithms as are adjoint sensitivities, they still are convenient 
for calculating gradients for objective function in the general form
\eqref{def:obj}. Specifically, the total derivative of the scalar function
$\Phi(X_N)$ can be computed with
\begin{equation}
\label{eqn:total_der}
\frac{d \Phi}{d p}(X_N) =  \frac{\partial \Phi}{\partial X}(X_N) \Sen_{N}  + \frac{\partial \Phi}{\partial p} (X_N).
\end{equation}
The total derivative of the integral term in \eqref{def:obj} (denoted by $q$ for simplicity) to parameters $p$ is given as
\begin{equation}
\label{eqn:total_der_p}
\frac{d q}{d p} = \int_{t_0}^{t_F} \left( \frac{\partial r}{\partial X}(t,X) \Sen + \frac{\partial r}{\partial p} (t,X) \right) \, dt.
\end{equation}
This integral together with $q$ is calculated automatically by PETSc with the
same timestepping algorithm and sequence of timesteps in the discrete
approaches for consistency, when users provide the necessary Jacobian callbacks
with {\tt TSSetCostIntegrand}. In addition, the forward approach is useful for
obtaining solution sensitivities often required by second-order adjoint
sensitivity analysis \cite{Barton_2005}.

\section{Handling discontinuities and events}\label{sec:tsevent}
One characteristic of applications, typically from the control systems
world, is the discontinuous nature of equations due to the presence of
various time- and state-based nonlinearities such as faults, limiters,
and hysterisis. Such discontinuities give rise to the following
conditionals\footnote{Equation (\ref{eq:disc}) shows one form for
  illustrative purposes. In general, the conditionals can include
  functions of the state variables instead of simple box constraints.}
introduced in the ODE or DAE equations:
\begin{equation}
	\begin{cases}
		x - x^+ = 0, & \text{if } x \ge x^+ \\
		x - x^-  = 0, & \text{if } x \le x^- \\
		\dot{x} = f(x), & \text{otherwise}.
	\end{cases}
	\label{eq:disc}
\end{equation}	
 
\begin{figure}[h]
  \centering
  \hspace*{-1cm}
  \includegraphics[scale=0.5]{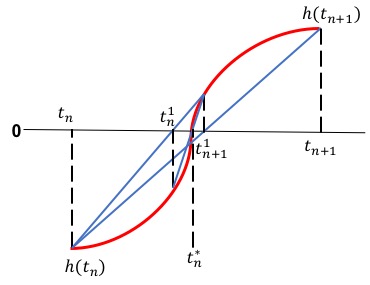}
  \caption{Detection and location of nonlinearities}
  \label{fig:event}
\end{figure}

PETSc supports the handling of such discontinuities through its
event-handling mechanism called {\em TSEvent}. Detecting and locating
such discontinuities are done by using an event handler or
root-finding method. A switching function $h(t,x) = 0$ is propagated
along with the equations. The mechanism of event detection and
location is illustrated in Fig. \ref{fig:event}.

The timestepper checks for the \emph{zero-crossing} of the event
function at every timestep. Specific directions of zero-crossing---positive
only, negative only, or both---can be provided. The
zero-crossing of an event is detected by the sign change of the event
function, namely, $sign(h(t_{n})) \ne sign(h(t_{n+1}))$. If this
condition is true, the event is said to be detected and the solution
rolled back to $t_n$. By using interpolation and successively shrinking
the time boundaries, the zero-crossing of the event function is
detected when its value is within a specified tolerance. At this time
instant, $t^*_n$ in Fig. \ref{fig:event}, the discontinuity is applied,
and an additional step is taken to synchronize with
$t_{n+1}$.
TSEvent also incorporates further
improvements to avoid duplicate steps (by utilizing the Illinois algorithm
\cite{Dowell1971}), and it speeds the detection of event zero-crossing
by using the Anderson--Bj\"orck method \cite{Galdino2011}. In the case of
multiple events detected during the same timestep, the event
detection mechanism uses the smallest interpolated timestep from the
list of events.

\begin{figure}[h]
  \centering
  \hspace*{-1cm}
  \includegraphics[scale=0.3]{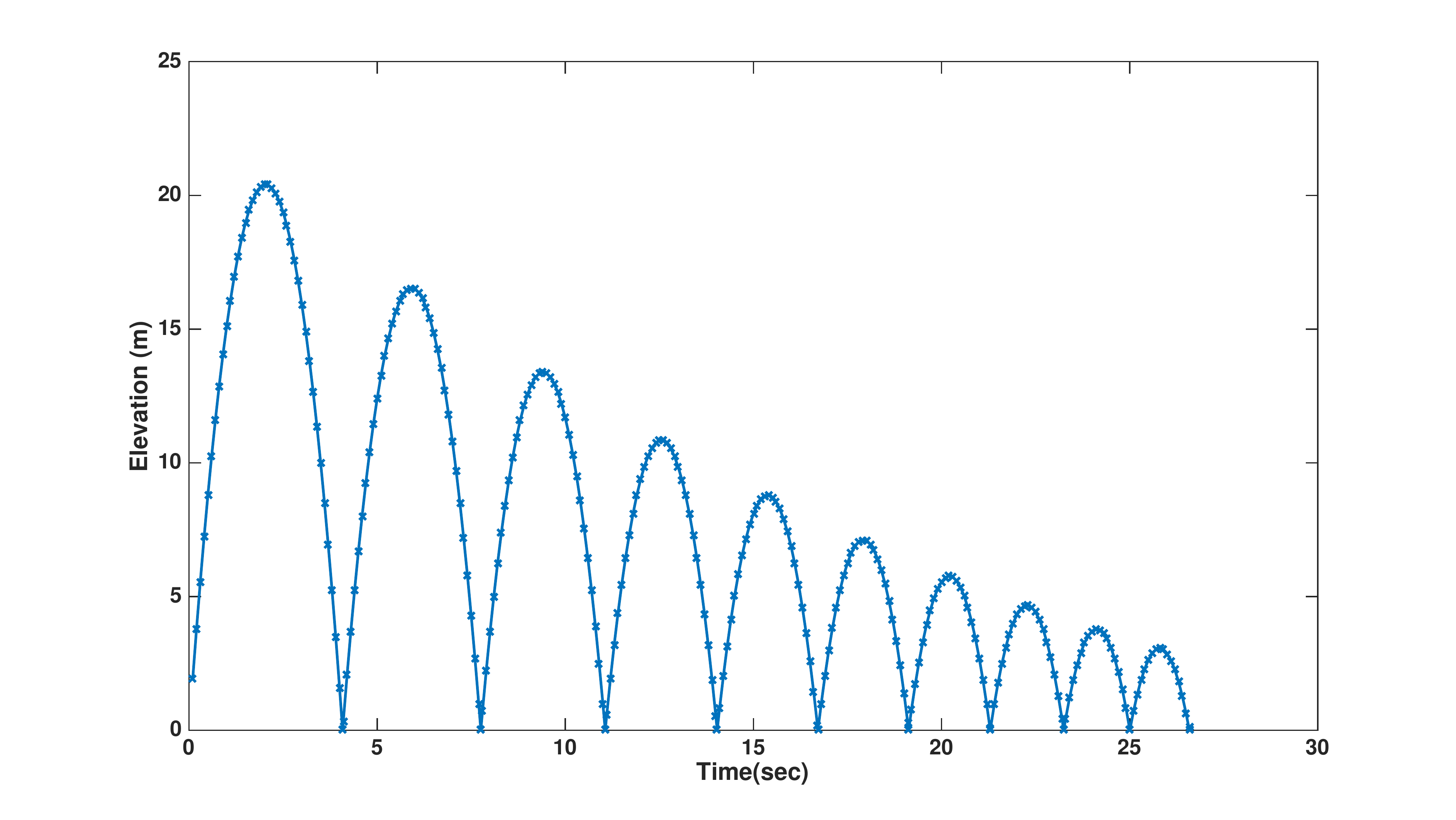}
  \caption{Bouncing ball example: The dynamics of the bouncing ball are described by the equations $\dot{u} = v$ and $\dot{v} = -9.8$. The ball velocity $v$ is attenuated by a factor of 0.9 every time it hits the ground $u = 0$.}
  \label{fig:event_example}
\end{figure}

Figure \ref{fig:event_example} presents a simple example illustrating
the usage of TSEvent for a bouncing ball. The event function is the
vertical position of the ball, $u$. When the ball hits the ground,
the sign change of the function is detected, and the discontinuity in the forcing function that
changes the sign of the velocity is applied, resulting in the ball
reversing direction.

Events can be set to TS through the application interface function {\em TSSetEventHandler()}, which has the following form:

 \begin{lstlisting}
TSSetEventHandler(TS ts,PetscInt nevents,PetscInt direction[],PetscBool terminate[],(*eventfun)(TS,PetscReal t,Vec X,PetscScalar h[],void *ctx),(*posteventfun)(TS ts,PetscInt nevents_det,PetscInt event_id[],PetscReal t,Vec X,PetscBool forwardsolve,void *ctx),void *ctx);
 \end{lstlisting}

 Here, \emph{nevents} is the number of local events to be located,
 \emph{direction[]} is an array of zero-crossing direction for each
 event, and \emph{terminate[]} array controls terminating TS
 timestepping after an event has been located. The event function
 $h(t,x)$ is set through the callback function \emph{*eventfun}; and,
 optionally, a post-event function \emph(*posteventfun) can be set
 that is called after an event or simultaneous multiple events are
 located. Specific actions following an event can be performed through
 the post-event function.

 For event functions having widely differing scales or range of values, finer control on locating the events can be provided through the {\em TSSetEventTolerances()} function:

 \begin{lstlisting}
TSSetEventTolerances(TS,PetscReal tol, PetscReal tols[]);
 \end{lstlisting}

 A single tolerance \emph{tol} can be used for all the events, or tolerances for each event can be set via the \emph{tols} array.
 
\section{Monitoring and Visualization}
Users of ODE solver packages often do not know much about even
the qualitative properties of the ODE they are solving; for example,
they may not even know whether it is stiff or which parts of the ODE are
stiff. 
To help users understand the qualitative properties of the solution,
PETSc/TS provides an extensible approach that allows
monitoring and visualizing the solution as well as solution properties,
such as maximum values of the solution or eigenvalues of the Jacobian.

Monitoring and visualization in PETSc are organized around the {\tt
  PetscViewer} object, which is an abstraction of ASCII and binary
files, as well as graphics APIs. Objects can be ``viewed'' with
varying levels of refinement based on the viewer used and options set
for the viewer. For example, {\tt MatView(A,viewer);} can display
minimalistically, with ASCII text, the size of the matrix and the
number of nonzeros or the entire matrix in binary format in a file or
as an image of the sparsity pattern depending on the viewer used. The
reduction of the object from its parallel representation is handled
automatically by PETSc. In addition to viewing the PETSc data objects
{\tt Vec} and {\tt Mat}, one can (in fact, doing so is desirable) view
the solver objects, for example, {\tt TS}. With an ASCII viewer it
prints information about the type of solver being used and all its
options; for binary viewers it saves the state of the object that can
be reloaded into memory with {\tt TSLoad()}; and for graphics viewers
it displays the relationship of the solver with the other solvers in
process, for example, that a SNES nonlinear solver object is embedded
in a TS object and that a KSP linear solver object is embedded in the
nonlinear solver object if Newton's method is being used.  Here we
display the output of a {\tt TSView()} on a particular ODE solver in
ASCII. The first part of the output summarizes the ODE integrator
information including the method used and its parameters. This is
followed by information about the linear  solver
(Rosenbrock-W methods solve only a linear system) which in this case is the direct solver
LU factorization.

{\small 
\begin{verbatim}
TS Object: 1 MPI processes
  type: rosw
  maximum steps=1000
  maximum time=20
  total number of nonlinear solver iterations=108
  total number of nonlinear solve failures=0
  total number of linear solver iterations=108
  total number of rejected steps=0
    Rosenbrock-W ra34pw2
    Abscissa of A       =  0.000000  0.871733  0.731580  1.000000 
    Abscissa of A+Gamma =  0.435867  0.871733  0.731580  1.000000 
  TSAdapt Object:   1 MPI processes
    type: basic
    number of candidates 1
      Basic: clip fastest decrease 0.1, fastest increase 10
      Basic: safety factor 0.9, extra factor after step rejection 0.5
    KSP Object:     1 MPI processes
      type: preonly
      maximum iterations=10000, initial guess is zero
      tolerances:  relative=1e-05, absolute=1e-50, divergence=10000
      left preconditioning
      using NONE norm type for convergence test
    PC Object:     1 MPI processes
      type: lu
        LU: out-of-place factorization
        tolerance for zero pivot 2.22045e-14
        matrix ordering: nd
        factor fill ratio given 5, needed 1
          Factored matrix follows:
            Mat Object:             1 MPI processes
              type: seqaij
              rows=3, cols=3
              package used to perform factorization: petsc
              total: nonzeros=9, allocated nonzeros=9
              total number of mallocs used during MatSetValues calls =0
                using I-node routines: found 1 nodes, limit used is 5
      linear system matrix = precond matrix:
      Mat Object:       1 MPI processes
        type: seqaij
        rows=3, cols=3
        total: nonzeros=9, allocated nonzeros=15
        total number of mallocs used during MatSetValues calls =0
          using I-node routines: found 1 nodes, limit used is 5
\end{verbatim}
}

Viewing of solver objects can usually be
controlled at runtime via the options database. For example, {\tt
  -ts\_view} produces ASCII output about the solver, whereas {\tt
  -ts\_view draw} produces a graphical display of the solver. 

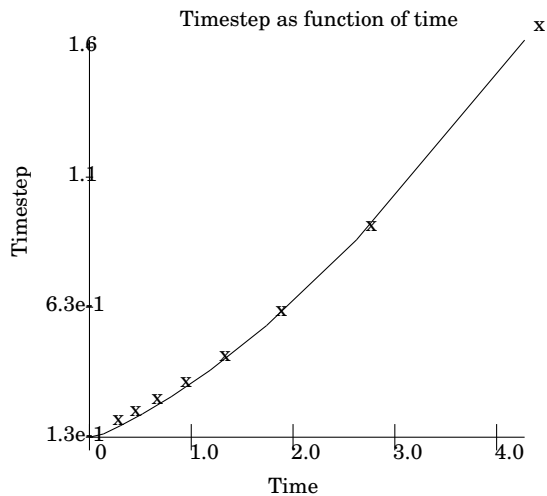
\begin{figure}[h]
\begin{frame}{}
\begin{changemargin}{-1cm}{0cm}
\begin{center}
\begin{tikzpicture}[scale = 7.00,font=\fontsize{8}{8}\selectfont]
\draw [black] (0.150413,0.169811) --(0.97686,0.169811);
\draw [black] (0.150413,0.169811) --(0.150413,0.924528);
\node [above right, black] at (0.304,0.924528) {Timestep as function of time};
\draw [black] (0.150413,0.144811) --(0.150413,0.194811);
\draw [black] (0.34372,0.144811) --(0.34372,0.194811);
\draw [black] (0.537027,0.144811) --(0.537027,0.194811);
\draw [black] (0.730333,0.144811) --(0.730333,0.194811);
\draw [black] (0.92364,0.144811) --(0.92364,0.194811);
\node [above right, black] at (0.143413,0.109811) {0};
\node [above right, black] at (0.32272,0.109811) {1.0};
\node [above right, black] at (0.516027,0.109811) {2.0};
\node [above right, black] at (0.709333,0.109811) {3.0};
\node [above right, black] at (0.90264,0.109811) {4.0};
\node [above right, black] at (0.472,0.0448113) {Time};
\draw [black] (0.143413,0.169811) --(0.157413,0.169811);
\draw [black] (0.143413,0.416949) --(0.157413,0.416949);
\draw [black] (0.143413,0.664086) --(0.157413,0.664086);
\draw [black] (0.143413,0.911223) --(0.157413,0.911223);
\node [above right, black] at (0.0496132,0.144811) {1.3e-1};
\node [above right, black] at (0.0496132,0.391949) {6.3e-1};
\node [above right, black] at (0.0916132,0.639086) {1.1};
\node [above right, black] at (0.0916132,0.886223) {1.6};
\node [rotate=90, black] at (0.021,0.599) {Timestep};
\draw [black] (0.150413,0.169811) --(0.176993,0.17599);
\node [above right, black] at (0.176993,0.17599) {x};
\draw [black] (0.176993,0.17599) --(0.209843,0.192024);
\node [above right, black] at (0.209843,0.192024) {x};
\draw [black] (0.209843,0.192024) --(0.251406,0.2143);
\node [above right, black] at (0.251406,0.2143) {x};
\draw [black] (0.251406,0.2143) --(0.305653,0.246734);
\node [above right, black] at (0.305653,0.246734) {x};
\draw [black] (0.305653,0.246734) --(0.379559,0.297001);
\node [above right, black] at (0.379559,0.297001) {x};
\draw [black] (0.379559,0.297001) --(0.486692,0.38196);
\node [above right, black] at (0.486692,0.38196) {x};
\draw [black] (0.486692,0.38196) --(0.657533,0.544857);
\node [above right, black] at (0.657533,0.544857) {x};
\draw [black] (0.657533,0.544857) --(0.97686,0.924528);
\node [above right, black] at (0.97686,0.924528) {x};
\end{tikzpicture}
\end{center}
\end{changemargin}
\end{frame}
\caption{Example of monitoring the adaptive timestep selected}
\label{dttime}
\end{figure}

In addition to static views of PETSc data and solver objects, we
provide numerous ways of dynamically viewing the solution and
properties of the solution, from within the program or via the
options database. This process is handled via ``monitor'' callback functions that can be attached to solver objects. For TS this is done with 
\begin{lstlisting}
TSMonitorSet(TS ts,(*monitor)(TS ts,PetscInt timestep,PetscReal time,Vec u ,void*mctx),void *mctx,(*mdestroy)(void**mctx));
\end{lstlisting}
The {\tt monitor()} function provided is called at the beginning and
at the end of
each timestep, and it can present the solution information in any way the
user likes. Various monitors may be set for the same
solver. PETSc provides a variety of default monitors that
\begin{itemize}
\item print the current timestep and time,
\item save the current solution to a binary or vtk file,
\item display the current solution by using a variety of graphical approaches using X windows or OpenGL, and
\item display the eigenvalues of the current operator, which is useful for understanding the stability of the scheme being used.
\end{itemize}
In addition, the {\tt monitor()} routines can compute and track
information over the lifetime of the simulation, for example, maximum
and minimum values of the solution or conserved quantities. The idea
is that rather than requiring users to modify the actual ODE
integrator code to track any property of the solution or solution
process, simple monitor routines are provided.  Many of these monitoring routines can be controlled from the command line; for example 
{\tt -ts\_monitor\_lg\_timestep} allows one to graphically monitor the changes in the adapted timestep as the computation proceeds as depicted in Fig. \ref{dttime}

The nonlinear and linear solvers also provide the same type of
flexible monitoring of the convergence process, with many available
default monitors allowing one to track how well the selected solvers
are working.

PETSc provides a simple but powerful API and options for gathering
performance information about the solution time, for example, time in
the linear solvers and time in the computation of the Jacobian. These
allow users to quickly focus in on the portions of the computation
that are the most time consuming and either select alternative
algorithms or further optimize the implementation. These are discussed in
the PETSc users manual \cite{petsc-user-ref}.

\section{Support for Specific Application Domains}
Many application areas have their own vocabulary and methodology for
describing their problem that are often distinct from the
language of ODEs. Although underlying their simulation is a set of ODEs
or DAEs, they never work directly with this form; rather, they express
their problems at a higher level of abstraction. Users can easily 
take advantage of these higher levels with the PETSc ODE and DAE
solvers, which allow the users to use their own natural language
for defining the problem and working with it. We already provide
this higher-level interface for two application areas: power systems
analysis and chemical reactions. We expect to do more in the future in
collaboration with application partners.

\subsection{Electrical Power Grid via DMNetwork}
Applications in an electrical power grid span a large range of temporal
and spatial scales that entail problems involving secure, stable, and
efficient planning and operation of the grid. A list of potential
applications suitable for PETSc usage is given in
\cite{AbhyankarHiPCNA11-1}. PETSc's timestepping library TS has been
used mainly for applications assessing the impacts of large
disturbances, such as short circuits and equipment outages, on the
stability of the grid. In such applications, called \emph{transient
  stability analysis} in electrical power grid parlance, the stability
of the grid is determined through a time-domain simulation of the
power grid differential-algebraic equations.
The differential
equations, $f(t,u,v)$, describe the dynamics of electromechanical
generators and motors, while the algebraic equations, $g(t,u,v)$, are
used for the electrical network comprising transmission lines,
transformers, and other connecting equipment.
\begin{align}
  \dot{u} =& f(t,u,v) \\
  0 =& g(t,u,v)
  \label{eq:epg_ts}
\end{align}
Different timestepping schemes, including adaptive
stepping and event handling, are compared in \cite{irep2017} for the
solution of transient stability problems. Rosenbrock schemes were found
to be optimal in terms of speedup and accuracy. In \cite{abhy2017}
the authors present experiments to achieve real-time simulation speed using
PETSc's timestepping and linear solvers. Results show that real-time
simulation speed was achieved on a fairly large electrical
grid. Similar real-time simulation speed results have been
reported in \cite{hipcna12}. \cite{jin2017} 
compare parallel transient stability algorithms using
MPI and OpenMP. GridPACK \cite{gridpack:homepage}, a software library
for developing parallel power grid applications that uses PETSc's core
solvers, is used in this work. Efficient calculation of sensitivities
of power grid dynamics trajectories to initial conditions using
a discrete adjoint scheme is described in \cite{Zhang_2017ab}.

DMNetwork \cite{siamns2014,mald2017} is a relatively new subclass of
PETSc's data management class DM that provides functionality for
efficiently managing and migrating data and topology for networks and
collections of networks. It handles the complex node-edge
relationships typically found in unstructured network problems; and it
provides simple abstractions to query the network topology and associate
physics with nodes/edges, acting as a middle layer between PETSc
solvers and the application physics. DMNetwork has been used for
several network applications, including electrical grids
\cite{asc2013}, water networks with over 1 billion unknowns \cite{mald2017}, and gas networks
\cite{jalving2017}.

\subsection{Chemical Reactions via TCHEM}

TCHEM \cite{safta2011tchem} is an open source implementation of many
of the reaction network chemistry capabilities of the commercial
ChemKin package \cite{chemkin-website}. TCHEM provides the code for
the ODE function evaluation and its Jacobian computation. It can read
ChemKin data files and construct the appropriate needed function
evaluations. PETSc provides an easy-to-use interface to TCHEM. In the
code listing below we demonstrate how the TCHEM function and Jacobian
routines can easily be wrapped and called from PETSc, thus merging
TCHEM's chemistry capabilities with PETSc's ODE integrators. The following code
segments demonstrate how to unwrap the PETSc data structure to
call TCHEM, which takes a raw array of numerical values on which to
apply the right-hand side function and compute the right-hand side
Jacobian and returns the Jacobian as a dense two-dimensional array
that is inserted back into the PETSc matrix.
\begin{lstlisting}
FormRHSFunction(TS ts,PetscReal t,Vec X,Vec F,void *ptr){
  User              user = (User)ptr;
  PetscScalar       *f;
  const PetscScalar *x;

  VecGetArrayRead(X,&x); VecGetArray(F,&f);
  PetscMemcpy(user->tchemwork,x,(user->Nspec+1)*sizeof(x[0]));
  user->tchemwork[0] *= user->Tini; /* Dimensionalize */
  TC_getSrc(user->tchemwork,user->Nspec+1,f);TC
  f[0] /= user->Tini;           /* Non-dimensionalize */
  VecRestoreArrayRead(X,&x); VecRestoreArray(F,&f);
}
FormRHSJacobian(TS ts,PetscReal t,Vec X,Mat Amat,Mat Pmat,void *ptr){
  User              user = (User)ptr;
  const PetscScalar *x;
  PetscInt          M = user->Nspec+1,i;

  VecGetArrayRead(X,&x);
  PetscMemcpy(user->tchemwork,x,(user->Nspec+1)*sizeof(x[0]));
  VecRestoreArrayRead(X,&x);
  user->tchemwork[0] *= user->Tini;  /* Dimensionalize temperature */
  TC_getJacTYN(user->tchemwork,user->Nspec,user->Jdense,1);
  for (i=0; i<M; i++) user->Jdense[i + 0*M] /= user->Tini; 
  for (i=0; i<M; i++) user->Jdense[0 + i*M] /= user->Tini; 
  for (i=0; i<M; i++) user->rows[i] = i;
  MatSetOption(Pmat,MAT_ROW_ORIENTED,PETSC_FALSE);
  MatSetOption(Pmat,MAT_IGNORE_ZERO_ENTRIES,PETSC_TRUE);
  MatZeroEntries(Pmat);
  MatSetValues(Pmat,M,user->rows,M,user->rows,user->Jdense,INSERT_VALUES);
  MatAssemblyBegin(Pmat,MAT_FINAL_ASSEMBLY);
  MatAssemblyEnd(Pmat,MAT_FINAL_ASSEMBLY);
}
\end{lstlisting}

\section{Conclusion}
PETSc provides a rich infrastructure for the efficient, scalable solution of ODEs and DAEs.
In this paper we have introduced the time integration component of
PETSc and provided examples of its usage.  We have listed and
described many of the integrators available and explained their basic
properties and usage.  In addition, we discussed the local and global
error control available with the integrators.  We also introduced the
capabilities for computing sensitivities (gradients) of functions of
the solutions to ODEs via forward and adjoint methods. We explained how events (discontinuities) in ODEs/DAEs may be handled and
the tools for monitoring and visualizing solutions and the solution
process. All the integrators are scalable and build on the basic mathematical libraries within PETSc.

{\bf Acknowledgment} We thank the many users of the TS for
  invaluable feedback on its capabilities and design. We thank
  Shashikant Aithal for his help in validating the TCHEM interfaces.
This material was based upon work supported by the U.S. Department of Energy, Office of
Science, Advanced Scientific Computing Research, under Contract
DE-AC02-06CH11357.

\bibliographystyle{ACM-Reference-Format-Journals}

\begin{thebibliography}{00}


\ifx \showCODEN    \undefined \def \showCODEN     #1{\unskip}     \fi
\ifx \showDOI      \undefined \def \showDOI       #1{{\tt DOI:}\penalty0{#1}\ }
  \fi
\ifx \showISBNx    \undefined \def \showISBNx     #1{\unskip}     \fi
\ifx \showISBNxiii \undefined \def \showISBNxiii  #1{\unskip}     \fi
\ifx \showISSN     \undefined \def \showISSN      #1{\unskip}     \fi
\ifx \showLCCN     \undefined \def \showLCCN      #1{\unskip}     \fi
\ifx \shownote     \undefined \def \shownote      #1{#1}          \fi
\ifx \showarticletitle \undefined \def \showarticletitle #1{#1}   \fi
\ifx \showURL      \undefined \def \showURL       #1{#1}          \fi

\bibitem[\protect\citeauthoryear{Abhyankar, Brown, Knepley, Meier, and
  Smith}{Abhyankar et~al\mbox{.}}{2014}]%
        {siamns2014}
{S. Abhyankar}, {J. Brown}, {M. Knepley}, {F. Meier}, {and} {B. Smith}. 2014.
\newblock \showarticletitle{Abstractions for expressing network problems in
  PETSc}. In {\em SIAM Conference on Network Science}.
\newblock


\bibitem[\protect\citeauthoryear{Abhyankar, Constantinescu, and
  Flueck}{Abhyankar et~al\mbox{.}}{2017a}]%
        {irep2017}
{S. Abhyankar}, {E. Constantinescu}, {and} {A. Flueck}. 2017a.
\newblock \showarticletitle{Variable-step multi-stage integration methods for
  fast and accurate power system dynamics simulation}. In {\em IREP'17}.
\newblock


\bibitem[\protect\citeauthoryear{Abhyankar, Constantinescu, Smith, Flueck, and
  Maldonado}{Abhyankar et~al\mbox{.}}{2017b}]%
        {abhy2017}
{S. Abhyankar}, {E. Constantinescu}, {B. Smith}, {A. Flueck}, {and} {D.
  Maldonado}. 2017b.
\newblock \showarticletitle{Parallel dynamics simulation using a
  {K}rylov-{S}chwarz linear solution scheme}.
\newblock {\em IEEE Transactions on Smart Grid\/}  {8} (2017), 1378--1386.
\newblock


\bibitem[\protect\citeauthoryear{Abhyankar and Flueck}{Abhyankar and
  Flueck}{2012}]%
        {hipcna12}
{S. Abhyankar} {and} {A. Flueck}. 2012.
\newblock \showarticletitle{Real-time power system dynamics simulation using a
  parallel block-{Jacobi} preconditioned {N}ewton-{GMRES} scheme}. In {\em
  HiPCNA-PG’12}.
\newblock


\bibitem[\protect\citeauthoryear{Abhyankar, Smith, Zhang, and Flueck}{Abhyankar
  et~al\mbox{.}}{2011}]%
        {AbhyankarHiPCNA11-1}
{Shrirang Abhyankar}, {Barry Smith}, {Hong Zhang}, {and} {A. Flueck}. 2011.
\newblock \showarticletitle{Using {PETSc} to Develop Scalable Applications for
  Next-Generation Power Grid}. In {\em Proceedings of the 1st International
  Workshop on High Performance Computing, Networking and Analytics for the
  Power Grid}. ACM.
\newblock
\showURL{%
\url{http://www.mcs.anl.gov/uploads/cels/papers/P1957-0911.pdf}}


\bibitem[\protect\citeauthoryear{Abhyankar, Smith, and
  Constantinescu}{Abhyankar et~al\mbox{.}}{2013}]%
        {asc2013}
{S. Abhyankar}, {B.~F. Smith}, {and} {E. Constantinescu}. 2013.
\newblock \showarticletitle{Evaluation of overlapping restricted additive
  Schwarz preconditioning for parallel solution of very large power flow
  problems}. In {\em HiPCNA-PG’13}.
\newblock


\bibitem[\protect\citeauthoryear{Ascher and Petzold}{Ascher and
  Petzold}{1998}]%
        {Ascher_B1998}
{U.M. Ascher} {and} {L.R. Petzold}. 1998.
\newblock {\em Computer Methods for Ordinary Differential Equations and
  Differential-Algebraic Equations}.
\newblock Society for Industrial Mathematics.
\newblock


\bibitem[\protect\citeauthoryear{Ascher, Ruuth, and Spiteri}{Ascher
  et~al\mbox{.}}{1997}]%
        {Ascher_1997}
{U.M. Ascher}, {S.J. Ruuth}, {and} {R.J. Spiteri}. 1997.
\newblock \showarticletitle{Implicit-explicit {R}unge-{K}utta methods for
  time-dependent partial differential equations}.
\newblock {\em Applied Numerical Mathematics\/}  {25} (1997), 151--167.
\newblock


\bibitem[\protect\citeauthoryear{Azyurt and Barton}{Azyurt and Barton}{2005}]%
        {Barton_2005}
{Derya~B. Azyurt} {and} {Paul~I. Barton}. 2005.
\newblock \showarticletitle{Cheap Second Order Directional Derivatives of Stiff
  ODE Embedded Functionals}.
\newblock {\em SIAM Journal on Scientific Computing\/} {26}, 5 (2005),
  1725--1743.
\newblock
\showDOI{%
\url{http://dx.doi.org/10.1137/030601582}}


\bibitem[\protect\citeauthoryear{Balay, Abhyankar, Adams, Brown, Brune,
  Buschelman, Dalcin, Eijkhout, Gropp, Kaushik, Knepley, May, McInnes, Mills,
  Munson, Rupp, Sanan, Smith, Zampini, Zhang, and Zhang}{Balay
  et~al\mbox{.}}{2018}]%
        {petsc-user-ref}
{Satish Balay}, {Shrirang Abhyankar}, {Mark~F. Adams}, {Jed Brown}, {Peter
  Brune}, {Kris Buschelman}, {Lisandro Dalcin}, {Victor Eijkhout}, {William~D.
  Gropp}, {Dinesh Kaushik}, {Matthew~G. Knepley}, {Dave~A. May}, {Lois~Curfman
  McInnes}, {Richard~Tran Mills}, {Todd Munson}, {Karl Rupp}, {Patrick Sanan},
  {Barry~F. Smith}, {Stefano Zampini}, {Hong Zhang}, {and} {Hong Zhang}. 2018.
\newblock {\em {PETS}c Users Manual}.
\newblock {T}echnical {R}eport ANL-95/11 - Revision 3.9. Argonne National
  Laboratory.
\newblock


\bibitem[\protect\citeauthoryear{Balay, Gropp, McInnes, and Smith}{Balay
  et~al\mbox{.}}{1997}]%
        {efficient}
{Satish Balay}, {William~D. Gropp}, {Lois~Curfman McInnes}, {and} {Barry~F.
  Smith}. 1997.
\newblock \showarticletitle{Efficient Management of Parallelism in Object
  Oriented Numerical Software Libraries}. In {\em Modern Software Tools in
  Scientific Computing}, {E.~Arge}, {A.~M. Bruaset}, {and} {H.~P. Langtangen}
  (Eds.). Birkhauser Press, 163--202.
\newblock


\bibitem[\protect\citeauthoryear{Boscarino, Pareschi, and Russo}{Boscarino
  et~al\mbox{.}}{2011}]%
        {Boscarino_TR2011}
{S. Boscarino}, {L. Pareschi}, {and} {G. Russo}. 2011.
\newblock Implicit-Explicit {R}unge-{K}utta schemes for hyperbolic systems and
  kinetic equations in the diffusion limit.  (2011).
\newblock
\newblock
\shownote{Arxiv preprint arXiv:1110.4375.}


\bibitem[\protect\citeauthoryear{Brenan, Campbell, Campbell, and
  Petzold}{Brenan et~al\mbox{.}}{1996}]%
        {Brenan_B1996}
{K.E. Brenan}, {S.L. Campbell}, {S.L.V. Campbell}, {and} {L.R. Petzold}. 1996.
\newblock {\em Numerical Solution of Initial-Value Problems in
  Differential-Algebraic Equations}.
\newblock Society for Industrial Mathematics.
\newblock


\bibitem[\protect\citeauthoryear{Butcher}{Butcher}{2008}]%
        {Butcher_B2008}
{J.C. Butcher}. 2008.
\newblock {\em Numerical Methods for Ordinary Differential Equations\/} (second
  ed.).
\newblock Wiley.
\newblock
\showISBNx{0470723351}


\bibitem[\protect\citeauthoryear{Butcher, Jackiewicz, and Wright}{Butcher
  et~al\mbox{.}}{2007}]%
        {Butcher_2007}
{J.C. Butcher}, {Z. Jackiewicz}, {and} {W.M. Wright}. 2007.
\newblock \showarticletitle{Error propagation of general linear methods for
  ordinary differential equations}.
\newblock {\em Journal of Complexity\/} {23}, 4-6 (2007), 560--580.
\newblock
\showISSN{0885-064X}
\showDOI{%
\url{http://dx.doi.org/10.1016/j.jco.2007.01.009}}


\bibitem[\protect\citeauthoryear{chemkin}{chemkin}{2017}]%
        {chemkin-website}
chemkin 2017.
\newblock ANSYS Chemkin-Pro Website.
\newblock \url{http://www.ansys.com/products/fluids/ansys-chemkin-pro}.
  (2017).
\newblock


\bibitem[\protect\citeauthoryear{Constantinescu}{Constantinescu}{2018}]%
        {Constantinescu_2018a}
{E.M. Constantinescu}. 2018.
\newblock \showarticletitle{Generalizing global error estimation for ordinary
  differential equations by using coupled time-stepping methods}.
\newblock {\it J. Comput. Appl. Math.} {332}, Supplement C (2018), 140--158.
\newblock
\showISSN{0377-0427}
\showDOI{%
\url{http://dx.doi.org/10.1016/j.cam.2017.05.012}}


\bibitem[\protect\citeauthoryear{Constantinescu and Sandu}{Constantinescu and
  Sandu}{2010}]%
        {Constantinescu_A2010a}
{E.M. Constantinescu} {and} {A. Sandu}. 2010.
\newblock \showarticletitle{Extrapolated implicit-explicit time stepping}.
\newblock {\em SIAM Journal on Scientific Computing\/} {31}, 6 (2010),
  4452--4477.
\newblock
\showDOI{%
\url{http://dx.doi.org/10.1137/080732833}}


\bibitem[\protect\citeauthoryear{Dalcin, Paz, Kler, and Cosimo}{Dalcin
  et~al\mbox{.}}{2011}]%
        {dalcin2011parallel}
{Lisandro~D Dalcin}, {Rodrigo~R Paz}, {Pablo~A Kler}, {and} {Alejandro Cosimo}.
  2011.
\newblock \showarticletitle{Parallel distributed computing using python}.
\newblock {\em Advances in Water Resources\/} {34}, 9 (2011), 1124--1139.
\newblock


\bibitem[\protect\citeauthoryear{Dowell and Jarratt}{Dowell and
  Jarratt}{1971}]%
        {Dowell1971}
{M. Dowell} {and} {P. Jarratt}. 1971.
\newblock \showarticletitle{A modified regula falsi method for computing the
  root of an equation}.
\newblock {\em BIT\/}  {11} (1971), 168--174.
\newblock


\bibitem[\protect\citeauthoryear{Galdino}{Galdino}{2011}]%
        {Galdino2011}
{Sergio Galdino}. 2011.
\newblock \showarticletitle{A family of regula falsi root-finding methods}. In
  {\em Proceedings of the 2011 World Congress on Engineering and Technology}.
\newblock


\bibitem[\protect\citeauthoryear{Gear}{Gear}{1971}]%
        {Gear_B1971}
{C.W. Gear}. 1971.
\newblock {\em Numerical Initial Value Problems in Ordinary Differential
  Equations}.
\newblock Prentice Hall PTR.
\newblock


\bibitem[\protect\citeauthoryear{Giraldo, Kelly, and Constantinescu}{Giraldo
  et~al\mbox{.}}{2013}]%
        {Giraldo_2013}
{F.X. Giraldo}, {J.F. Kelly}, {and} {E.M. Constantinescu}. 2013.
\newblock \showarticletitle{Implicit-explicit formulations of a
  three-dimensional nonhydrostatic unified model of the atmosphere {(NUMA)}}.
\newblock {\em SIAM Journal on Scientific Computing\/} {35}, 5 (2013),
  B1162--B1194.
\newblock
\showDOI{%
\url{http://dx.doi.org/10.1137/120876034}}


\bibitem[\protect\citeauthoryear{Griewank and Walther}{Griewank and
  Walther}{2000}]%
        {griewank2000algorithm}
{Andreas Griewank} {and} {Andrea Walther}. 2000.
\newblock \showarticletitle{Algorithm 799: revolve: an implementation of
  checkpointing for the reverse or adjoint mode of computational
  differentiation}.
\newblock {\it ACM Trans. Math. Software} {26}, 1 (2000), 19--45.
\newblock


\bibitem[\protect\citeauthoryear{Hairer, N{\o}rsett, and Wanner}{Hairer
  et~al\mbox{.}}{2008}]%
        {Hairer_B2008_I}
{E. Hairer}, {S.P. N{\o}rsett}, {and} {G. Wanner}. 2008.
\newblock {\em Solving Ordinary Differential Equations {I}: {N}onstiff
  Problems}.
\newblock Springer.
\newblock
\showISBNx{978-3-540-78862-1}
\showDOI{%
\url{http://dx.doi.org/10.1007/978-3-540-78862-1}}


\bibitem[\protect\citeauthoryear{Hairer and Wanner}{Hairer and Wanner}{2002}]%
        {Hairer_B2002_II}
{E. Hairer} {and} {G. Wanner}. 2002.
\newblock {\em Solving Ordinary Differential Equations {II}: Stiff and
  Differential-Algebraic Problems}.
\newblock Springer.
\newblock
\showISBNx{3-540-60452-9}


\bibitem[\protect\citeauthoryear{Heroux, Bartlett, Hoekstra, Hu, Kolda,
  Lehoucq, Long, Pawlowski, Phipps, Salinger, Thornquist, Tuminaro,
  Willenbring, and Williams}{Heroux et~al\mbox{.}}{2003}]%
        {Trilinos}
{M. Heroux}, {R. Bartlett}, {V.H.R. Hoekstra}, {J. Hu}, {T. Kolda}, {R.
  Lehoucq}, {K. Long}, {R. Pawlowski}, {E. Phipps}, {A. Salinger}, {H.
  Thornquist}, {R. Tuminaro}, {J. Willenbring}, {and} {A. Williams}. 2003.
\newblock {\em An Overview of {T}rilinos}.
\newblock {T}echnical {R}eport SAND2003-2927. Sandia National Laboratories.
\newblock


\bibitem[\protect\citeauthoryear{Hindmarsh, Brown, Grant, Lee, Serban,
  Shumaker, and Woodward}{Hindmarsh et~al\mbox{.}}{2005}]%
        {Sundials}
{A.C. Hindmarsh}, {P.N. Brown}, {K.E. Grant}, {S.L. Lee}, {R. Serban}, {D.E.
  Shumaker}, {and} {C.S. Woodward}. 2005.
\newblock \showarticletitle{{SUNDIALS}: Suite of nonlinear and
  differential/algebraic equation solvers}.
\newblock {\em ACM Transactions on Mathematical Software (TOMS)\/} {31}, 3
  (2005), 363--396.
\newblock


\bibitem[\protect\citeauthoryear{Horn}{Horn}{1983}]%
        {Horn_1983}
{M.K. Horn}. 1983.
\newblock \showarticletitle{Fourth-and fifth-order, scaled Runge-Kutta
  algorithms for treating dense output}.
\newblock {\it SIAM J. Numer. Anal.} {20}, 3 (1983), 558--568.
\newblock


\bibitem[\protect\citeauthoryear{Jalving, Abhyankar, Kim, Herald, and
  Zavala}{Jalving et~al\mbox{.}}{2017}]%
        {jalving2017}
{J. Jalving}, {S. Abhyankar}, {K. Kim}, {M. Herald}, {and} {V. Zavala}. 2017.
\newblock \showarticletitle{A graph-based computational framework for
  simulation and optimization of coupled infrastructure networks}.
\newblock {\em IET Generation, Transmission, and Distribution\/}  {11} (2017),
  3163--3176.
\newblock


\bibitem[\protect\citeauthoryear{Jansen, Whiting, and Hulbert}{Jansen
  et~al\mbox{.}}{2000}]%
        {Jansen_2000}
{K.E. Jansen}, {C.H. Whiting}, {and} {G.M. Hulbert}. 2000.
\newblock \showarticletitle{A generalized-$\alpha$ method for integrating the
  filtered {N}avier--{S}tokes equations with a stabilized finite element
  method}.
\newblock {\em Computer Methods in Applied Mechanics and Engineering\/} {190},
  3 (2000), 305--319.
\newblock


\bibitem[\protect\citeauthoryear{Jin, Huang, Diao, Wu, and Chen}{Jin
  et~al\mbox{.}}{2017}]%
        {jin2017}
{S. Jin}, {Z. Huang}, {R. Diao}, {D. Wu}, {and} {Y. Chen}. 2017.
\newblock \showarticletitle{Comparative implementation of high performance
  computing for power system dynamic simulations}.
\newblock {\em IEEE Transactions on Smart Grid\/}  {8} (2017), 1387--1395.
\newblock


\bibitem[\protect\citeauthoryear{Kaps, Poon, and Bui}{Kaps
  et~al\mbox{.}}{1985}]%
        {Kaps_1985}
{P. Kaps}, {S.W.H. Poon}, {and} {T.D. Bui}. 1985.
\newblock \showarticletitle{Rosenbrock methods for stiff {ODE}s: {A} comparison
  of {R}ichardson extrapolation and embedding technique}.
\newblock {\em Computing\/} {34}, 1 (1985), 17--40.
\newblock


\bibitem[\protect\citeauthoryear{Kennedy and Carpenter}{Kennedy and
  Carpenter}{2003}]%
        {Kennedy_2003}
{C.A. Kennedy} {and} {M.H. Carpenter}. 2003.
\newblock \showarticletitle{Additive {R}unge-{K}utta schemes for
  convection-diffusion-reaction equations}.
\newblock {\em Appl. Numer. Math.\/} {44}, 1-2 (2003), 139--181.
\newblock
\showDOI{%
\url{http://dx.doi.org/10.1016/S0168-9274(02)00138-1}}


\bibitem[\protect\citeauthoryear{Ketcheson}{Ketcheson}{2008}]%
        {Ketcheson_2008}
{D.I. Ketcheson}. 2008.
\newblock \showarticletitle{Highly Efficient Strong Stability-Preserving
  {R}unge--{K}utta Methods with Low-Storage Implementations}.
\newblock {\em SIAM Journal on Scientific Computing\/} {30}, 4 (2008),
  2113--2136.
\newblock
\showDOI{%
\url{http://dx.doi.org/10.1137/07070485X}}


\bibitem[\protect\citeauthoryear{Maldonado, Abhyankar, Smith, and
  Zhang}{Maldonado et~al\mbox{.}}{2017}]%
        {mald2017}
{D. Maldonado}, {S. Abhyankar}, {B. Smith}, {and} {H. Zhang}. 2017.
\newblock Scalable multiphysics network simulation using {PETSc DMNetwork}.
\newblock \url{http://www.mcs.anl.gov/papers/P7065-0617.pdf}.   (2017).
\newblock


\bibitem[\protect\citeauthoryear{Marin, Constantinescu, and Smith}{Marin
  et~al\mbox{.}}{2017}]%
        {pdeconstrainedspectraladjoints}
{Oana Marin}, {Emil Constantinescu}, {and} {Barry Smith}. 2017.
\newblock {\em PDE constrained optimization, error estimation and control with
  spectral elements using PETSc and TAO}.
\newblock Preprint ANL/MCS-P9031-1117. ANL.
\newblock


\bibitem[\protect\citeauthoryear{MATLAB}{MATLAB}{2014}]%
        {MATLAB}
{MATLAB}. 2014.
\newblock {\em version 8.1.0 (R2013a)}.
\newblock The MathWorks Inc., Natick, Massachusetts.
\newblock


\bibitem[\protect\citeauthoryear{NAG}{NAG}{2018}]%
        {NAG}
{NAG}. 2018.
\newblock The NAG Library, The Numerical Algorithms Group (NAG).
\newblock Oxford, United Kingdom.   (2018).
\newblock
\showURL{%
\url{www.nag.com}}


\bibitem[\protect\citeauthoryear{{Palmer et al.}}{{Palmer et al.}}{2018}]%
        {gridpack:homepage}
{B. {Palmer et al.}} 2018.
\newblock {GridPACK {W}eb page}.
\newblock \url{https://www.gridpack.org/wiki/index.php/Main_Page}.   (2018).
\newblock


\bibitem[\protect\citeauthoryear{Pareschi and Russo}{Pareschi and
  Russo}{2005}]%
        {Pareschi_2005}
{L. Pareschi} {and} {G. Russo}. 2005.
\newblock \showarticletitle{Implicit-Explicit {R}unge-{K}utta Schemes and
  Applications to Hyperbolic Systems with Relaxation}.
\newblock {\em Journal of Scientific Computing\/} {25}, 1 (2005), 129--155.
\newblock


\bibitem[\protect\citeauthoryear{Petzold}{Petzold}{1992}]%
        {DASSL}
{L.R. Petzold}. 1992.
\newblock \showarticletitle{DASSL, Solution of Differential Algebraic
  Equation}.
\newblock  (1992).
\newblock


\bibitem[\protect\citeauthoryear{Rang and Angermann}{Rang and
  Angermann}{2005}]%
        {Rang_2005}
{J. Rang} {and} {L. Angermann}. 2005.
\newblock \showarticletitle{New {R}osenbrock {W}-methods of order 3 for partial
  differential algebraic equations of index 1}.
\newblock {\em BIT Numerical Mathematics\/} {45}, 4 (2005), 761--787.
\newblock


\bibitem[\protect\citeauthoryear{Safta, Najm, and Knio}{Safta
  et~al\mbox{.}}{2011}]%
        {safta2011tchem}
{Cosmin Safta}, {Habib~N Najm}, {and} {Omar Knio}. 2011.
\newblock \showarticletitle{TChem-a software toolkit for the analysis of
  complex kinetic models}.
\newblock {\em Sandia Report, SAND2011-3282\/} (2011).
\newblock


\bibitem[\protect\citeauthoryear{Sandu, Verwer, Blom, Spee, Carmichael, and
  Potra}{Sandu et~al\mbox{.}}{1997}]%
        {Sandu_1997}
{A. Sandu}, {J.G. Verwer}, {J.G. Blom}, {E.J. Spee}, {G.R. Carmichael}, {and}
  {F.A. Potra}. 1997.
\newblock \showarticletitle{Benchmarking stiff ode solvers for atmospheric
  chemistry problems {II}: {R}osenbrock solvers}.
\newblock {\em Atmospheric Environment\/} {31}, 20 (1997), 3459--3472.
\newblock


\bibitem[\protect\citeauthoryear{S{\"o}derlind}{S{\"o}derlind}{2003}]%
        {Soderlind_2003}
{G. S{\"o}derlind}. 2003.
\newblock \showarticletitle{Digital filters in adaptive time-stepping}.
\newblock {\em ACM Transactions on Mathematical Software (TOMS)\/} {29}, 1
  (2003), 1--26.
\newblock


\bibitem[\protect\citeauthoryear{S{\"o}derlind and Wang}{S{\"o}derlind and
  Wang}{2006}]%
        {Soderlind_2006}
{G. S{\"o}derlind} {and} {L. Wang}. 2006.
\newblock \showarticletitle{Adaptive time-stepping and computational
  stability}.
\newblock {\it J. Comput. Appl. Math.} {185}, 2 (2006), 225--243.
\newblock


\bibitem[\protect\citeauthoryear{Zhang, Abhyankar, Constantinescu, and
  Anitescu}{Zhang et~al\mbox{.}}{2017}]%
        {Zhang_2017ab}
{H. Zhang}, {S.S. Abhyankar}, {E.M. Constantinescu}, {and} {Mihai Anitescu}.
  2017.
\newblock \showarticletitle{Discrete adjoint sensitivity analysis of hybrid
  dynamical systems with switching}.
\newblock {\em IEEE Transactions on Circuits and Systems I: Regular Papers\/}
  {64}, 5 (2017), 1247--1259.
\newblock
\showDOI{%
\url{http://dx.doi.org/10.1109/TCSI.2017.2651683}}


\bibitem[\protect\citeauthoryear{Zhong}{Zhong}{1996}]%
        {Zhong_1996}
{X. Zhong}. 1996.
\newblock \showarticletitle{Additive semi-implicit {R}unge-{K}utta methods for
  computing high speed nonequilibrium reactive flows}.
\newblock {\it J. Comput. Phys.}  {128} (1996), 19--31.
\newblock


\end{thebibliography}


\newpage
{\bf Disclaimer.} The submitted manuscript has been created by UChicago Argonne, LLC,
Operator of Argonne National Laboratory (``Argonne'').
Argonne, a U.S. Department of Energy Office of Science laboratory, is
operated under Contract No. DE-AC02-06CH11357. The U.S. Government
retains for itself, and others acting on its behalf, a paid-up
nonexclusive, irrevocable worldwide license in said article to reproduce,
prepare derivative works, distribute copies to the public, and perform
publicly and display publicly, by or on behalf of the Government.
The Department of Energy will provide public access to these results of federally sponsored research in accordance with the DOE Public Access Plan. http://energy.gov/downloads/doe-public-access-plan.

\end{document}